\documentclass[10pt,draft,twoside]{amsart}
\usepackage{amssymb}
\usepackage{latexsym}
\usepackage{amsfonts}
\usepackage{amsmath}
\usepackage{fancyheadings}
\usepackage{bookman}

\newcommand{\K}{\mathcal K}
\newcommand{\T}{\mathcal T}

\newcommand{\E}{\mathcal E}

\newcommand{\cd}[2]{{\sf CD}_{\rm #1}(#2)}

\newtheorem{theorem}{Theorem}[section]

\newtheorem{lemma}[theorem]{Lemma}
\newtheorem{corollary}[theorem]{Corollary}

\theoremstyle{definition}

\newtheorem{example}[theorem]{Example}

\newcommand{\alt}[1]{{\sf A}_{#1}}
\newcommand{\Alt}[1]{{\sf Alt}\,#1}
\newcommand{\mat}[1]{{\sf M}_{#1}}
\newcommand{\sy}[1]{{\sf S}_{#1}}

\renewcommand{\sp}[2]{{\sf Sp}_{#1}(#2)}
\newcommand{\pomegap}[2]{{\sf P}\Omega^+_{#1}(#2)}

\newcommand{\sym}[1]{{\sf Sym}\,#1}
\newcommand{\supp}{\operatorname{Supp}}

\renewcommand{\wr}{\,{\sf wr}\,}

\newcommand{\dih}[1]{{\sf D}_{#1}}

\newcommand{\aut}[1]{{\sf Aut}\,{#1}}

\newcommand{\cent}[2]{{\mathbb C}_{#1}(#2)}
\newcommand{\norm}[2]{{\mathbb N}_{#1}\left(#2\right)}

\newcommand{\psl}[2]{\mbox{\sf PSL}_{#1}(#2)}

\renewcommand{\leq}{\leqslant}
\renewcommand{\geq}{\geqslant}

\renewcommand{\L}{\mathcal L}

\begin{document}

\title{Innately transitive subgroups of wreath products in product action}
\author{Robert W. Baddeley, Cheryl E. Praeger and Csaba Schneider}
\address[Baddeley]{32 Arbury Road\\Cambridge CB4 2JE, UK}
\address[Praeger and Schneider]{School of Mathematics and Statistics\\
The University of Western Australia\\
35 Stirling Highway 6009 Crawley\\
Western Australia}
\curraddr[Schneider]{Informatics Laboratory\\ Computer and Automation Research Institute of the Hungarian Academy of Sciences\\ L\'agym\'anyosi u. 11.\\1111 Budapest, Hungary}
\email{robert.baddeley@ntworld.com, praeger@maths.uwa.edu.au,\hfill\break
csaba@maths.uwa.edu.au\protect{\newline} {\it WWW:}
www.maths.uwa.edu.au/$\sim$praeger, www.maths.uwa.edu.au/$\sim$csaba}

\begin{abstract}
A permutation group is innately transitive if it has a transitive
minimal normal subgroup, which is referred to as a plinth. 
We study the class of finite, innately
transitive permutation groups that can be embedded into wreath
products in product action. This investigation is carried out by
observing that such a wreath product preserves a natural Cartesian
decomposition 
of the underlying set. Previously we classified the possible
embeddings in the case where the plinth is simple. Here we extend that
classification and identify several different types of Cartesian
decompositions that can be preserved by an innately transitive group
with a non-abelian plinth. These different types of decompositions
lead to different types of embeddings of the acting group 
into wreath products in product action. We also obtain a full
characterisation of embeddings of innately transitive groups with
diagonal type into such wreath products.
\end{abstract}

\thanks{{\it Date:} draft typeset \today
The authors acknowledge the support of an Australian Research Council grant. 
We are very grateful to Laci Kov\'acs for explaining the origins of some 
of the ideas that appear in this paper.}

\maketitle

\section{Introduction}

This paper forms an important part of our program to describe innately
transitive subgroups of wreath products in product action (see Section~\ref{cssect} for the definition of
wreath products and their product actions). A
permutation group is {\em innately transitive} if it has a transitive
minimal normal subgroup; and such a 
subgroup is called a {\em plinth}. A
permutation group is said to be {\em quasiprimitive} if all its minimal normal
subgroups are transitive. 
There are various
characterisations of innately transitive and quasiprimitive
groups: in~\cite{bp} 5~principal types of finite innately transitive groups
were identified, and~\cite{bad:quasi} listed 8 types of finite
quasiprimitive permutation groups based on the classification
obtained by~\cite{prae:quasi}.
We use the types of~\cite{bp}
and~\cite{bad:quasi} in this paper.

In an earlier paper~\cite{recog} we described those finite,
innately transitive groups with a
simple plinth that can be embedded into a wreath product in product
action. The main results of this paper extend that classification to
embeddings of 
innately transitive subgroups with a
non-abelian plinth in such wreath products. 
Let $H$ be a quasiprimitive,  almost simple permutation group acting
on a set $\Gamma$. That is, $H$ has a unique minimal normal subgroup
$U$, and $U$ is a non-abelian simple group.  Set $W=H\wr\sy \ell$ for
some $\ell\geq 2$, and
consider $W$ as a permutation
group in product action acting on $\Gamma^\ell$. Let $N=U_1\times\cdots\times
U_\ell$ be the unique minimal normal subgroup of $W$; note that $N\cong U^\ell$.  
Assume that $G$ is an innately transitive subgroup of $W$ with a
non-abelian plinth $M=T_1\times\cdots\times T_k$ where $T_1,\ldots,T_k$ are
finite, non-abelian simple groups all isomorphic to a group $T$. 

The following two theorems are the main results of this paper.

\begin{theorem}\label{main}
If $G$, $W$, $M$, and $N$ are as above, then $M\leq N$. Further, if $G$ projects onto a transitive
subgroup of $\sy\ell$, then exactly one of
the following holds.
\begin{enumerate}
\item[(a)] $k=\ell$; the $T_i$ and the $U_i$ can be indexed
so that $T_1\leq U_1,\ldots,T_k\leq U_k$.
\item[(b)] $\ell=2k$; $T$ and $U$ are as in Table~$\ref{maintable}$; the $T_i$ and the $U_i$ can be indexed
so that $T_1\leq U_1\times U_2,T_2\leq U_3\times U_4,\ldots,T_k\leq
U_{2k-1}\times U_{2k}$. 
\item[(c)] None of the cases~(a)--(b) holds and $U=\Alt \Gamma$.
\end{enumerate}
\end{theorem}

\begin{center}
\begin{table}
$$
\begin{array}{|c|c|c|}
\hline
& T & U \\
\hline\hline
1& \alt 6 & \alt 6 \\
\hline
2 & \mat{12} & \mat{12} \\
\hline
3 & \mat{12} & \alt{12} \\
\hline
4 & \pomegap 8q & \pomegap 8q \\
\hline
5 & \pomegap 8q & \alt n\mbox{ where } n=|\pomegap 8q:\Omega_7(q)|\\
\hline
6 & \pomegap 82 & \sp 82 \\\hline
7 & \sp 4{2^a},\ a\geq 2 & \sp{4b}{2^{a/b}} \mbox{ where }\ b\mid a \\
\hline
8 & \sp 4{2^a},\ a\geq 2 & \alt n \mbox{ where } n=|\sp 4{2^a}:\sp 2{2^{2a}}\cdot 2|\\
\hline
\end{array}
$$
\caption{Table for Theorem~\ref{main}}\label{maintable}
\end{table}
\end{center}

An even stronger result can be obtained if $G$ has diagonal
type, which is defined as follows. Let $\sigma_i:M\rightarrow T_i$ be
the natural projection map, for $i\in\{1,\ldots,k\}$. An innately transitive group $G$
has {\em diagonal type} if $\sigma_i(M_\omega)=T_i$ for all $\omega\in\Gamma^\ell$ and 
$i\in\{1,\ldots,k\}$. It
follows from Scott's Lemma~\ref{scott} that in this case a point
stabiliser $M_\omega$ is isomorphic to $T^s$ for some $s\geq 1$. If
$s=1$ then we say that $G$ has {\em simple diagonal type}, otherwise
$G$ has {\em compound diagonal type}. It was shown
in~\cite[Proposition~5.5]{bp} that an innately transitive group of
diagonal type contains a unique minimal normal subgroup, and hence is quasiprimitive.

\begin{theorem}\label{diagonalembed}
If $W$, $G$, $N$, $M$ are as above and $G$ has diagonal type then the
following all hold.
\begin{enumerate}
\item[(a)] $G$ projects onto a transitive subgroup of $\sy\ell$ and
$G$ is a quasiprimitive group of compound diagonal type;
\item[(b)] $U= \Alt\Gamma$;
\item[(c)] $M\leq N$; $k=m\ell$ for some $m\geq 2$ 
and the $U_i$ and the $T_i$ can be indexed so that
$T_1\times\cdots\times T_m\leq
U_1$, $T_{m+1}\times\cdots\times T_{2m}\leq
U_2,\ldots,T_{(\ell-1)m+1}\times\cdots\times T_{\ell m}\leq U_\ell$.
\end{enumerate}
\end{theorem}

Theorem~\ref{diagonalembed}
implies the following important corollary, thus
proving~\cite[Theorem~4.7(3)]{bad:quasi}. 

\begin{corollary}\label{sdcor}
An innately transitive permutation group with simple diagonal type can
never be a subgroup of  a wreath product $\sym\Gamma\wr\sy\ell$ in
product action with
$|\Gamma|\geq 2$ and $\ell\geq 2$. 
\end{corollary}

Note that cases (a) and (b) of Theorem~\ref{main} give detailed
information about the embedding $G\leq W$, while case~(c) contains a
rich variety of examples which will be
investigated further in Section~\ref{final}. The
inclusions in case~(a) are quite common. Take, for instance, finite
simple permutation groups $T,\ U\leq \sym\Gamma$ such that
$T\leq U$. If $K$ is a transitive subgroup of $\sy\ell$ then clearly
$T\wr K\leq U\wr\sy\ell\leq\sym{\Gamma^\ell}$, and this inclusion is
as in Theorem~\ref{main}(a). 

Embeddings belonging
to Theorem~\ref{main}(b) are not a great deal more mysterious, as illustrated
by the following construction. 
By~\cite{lps:maxsub}, 
the group $\aut{\alt 6}\cong{\sf P}\Gamma{\sf L}_2(9)$, acting on a set
of size 36, can be embedded into
$\sy 6\wr\sy 2$ such that $\aut{\alt 6}$ projects onto $\sy 2$ via the
natural projection map $\sy 6\wr\sy 2\rightarrow\sy 2$.
Hence if $K$ is a transitive subgroup of $\sy k$ for some $k\geq
2$ then $(\aut{\alt 6})\wr
K$ is a subgroup of $(\sy 6\wr\sy 2)\wr \sy k$, and, in turn, of $\sy 6\wr\sy{2k}$ in product action. It is easy to see that this
inclusion $(\aut{\alt 6})\wr
K\leq \sy 6\wr\sy {2k}$ is as in Theorem~\ref{main}(b). 

The easiest examples in case~(c) are constructed by taking a simple
permutation group $T\leq\sym\Gamma$ and a transitive group $K\leq \sy
m$. Then $T\wr K$ is a subgroup of $\sym{\Gamma}\wr\sy m$ in product
action on $\Gamma^m$. Moreover if $L\leq\sy\ell$ is a transitive permutation group
then clearly $(T\wr K)\wr L\leq
\sym{\Gamma^m}\wr\sy\ell\leq\sym{\Gamma^{m\ell}}$. The inclusion
$(T\wr K)\wr L\leq \sym{\Gamma^m}\wr\sy\ell$ is as in
Theorem~\ref{main}(c). However, as we will see in Section~\ref{exsec},
the most intriguing examples in case~(c) cannot be obtained by this
simple construction.

Note that, in Theorem~\ref{diagonalembed}, the fact that $G$ induces a
transitive subgroup of $\sy\ell$ is a consequence of the
hypotheses. 
An easy example of the inclusions described by
Theorem~\ref{diagonalembed} can be constructed as follows. Let $K$ be
an innately transitive permutation group 
of simple diagonal type acting on $\Gamma$ and let $L$ be its unique minimal normal
subgroup. 
Let $T$ denote the isomorphism type of a simple direct factor of
$M$. Then for $\gamma\in\Gamma$, the point stabiliser $L_\gamma$ is
isomorphic to $T$. If $Q$ is a transitive subgroup of $\sy\ell$, then
$K \wr Q$ is an innately transitive group of compound diagonal type
acting on $\Gamma^\ell$, and clearly $K\wr Q\leq \sym\Gamma\wr\sy\ell$. This
inclusion gives another example for Theorem~\ref{main}(c).

In Section~\ref{exsec} we give three further examples of inclusions of
innately transitive groups into wreath products in product action to
illustrate the diversity and the beauty of the embeddings belonging to
Theorem~\ref{main}(c). Section~\ref{cssect} contains a brief account of
Cartesian decompositions of sets and Cartesian systems of subgroups
that were introduced in~\cite{recog} and are fundamental to our
approach. 

The concept of Cartesian decompositions 
can be viewed as a straightforward generalisation of
the concept of systems of product imprimitivity. Such systems 
were first defined by L.\ G.\ Kov\'acs in his 
lecture at the Group Theory Conference at Oberwolfach in May 1987. The ideas 
presented in this Oberwolfach lecture played a crucial r\^ole in his paper~\cite{kov:decomp}.  
In the terminology used by Kov\'acs, a system of product imprimitivity
for a group $G$ is a special type of Cartesian decomposition preserved by $G$, namely it is one on which $G$ acts transitively.

We draw attention to Theorem~\ref{disstrips}, which is  of crucial
importance in analysing the structure of Cartesian systems in $M$ in which
the elements involve diagonal subgroups (which are called strips
in our terminology introduced in Section~\ref{charsec}) isomorphic to $T$, as in
Example~\ref{ex1s}.
This result may be
compared with Scott's Lemma~\ref{scott} on subdirect
subgroups in characteristically simple groups that played a key r\^ole
in proving the O'Nan-Scott Theorem for primitive groups.

Section~\ref{4} contains our major result, the 6-Class Theorem, about
transitive Cartesian decompositions preserved by a finite innately
transitive group. The examples in the introduction and in
Section~\ref{exsec} show that the classification in Theorem~\ref{main}
can naturally be refined using the language of Cartesian systems. 
If $G\leq\sym\Gamma\wr\sy\ell$ and $G$ induces a transitive subgroup
of $\sy\ell$, then we say that the corresponding Cartesian
decomposition is transitive.
We find that for
an innately transitive group $G$ with a non-abelian plinth, each of
the transitive $G$-invariant Cartesian decompositions of the underlying set
belongs to exactly one of the 6 classes identified by the 6-Class
Theorem (Theorem~\ref{5class}). The class of a particular Cartesian decomposition can be
determined from the abstract group structure of the subgroups in the
corresponding Cartesian system. Cartesian decompositions in different
classes lead to different types of embeddings of the group $G$ into
wreath products in product action. This is  strong evidence to support
our belief that the theory
of Cartesian decompositions and Cartesian systems is the appropriate
framework for studying these embeddings.

If $G$ is an innately transitive permutation group with a non-abelian
plinth, each transitive $G$-invariant Cartesian decomposition 
gives rise to a certain interesting combinatorial
structure preserved by the group or its stabiliser subgroups. These
combinatorial structures may be partitions (see Theorem~\ref{c1}), bipartite graphs,
generalised graphs, depending on the class of the decomposition. In
some classes it is also possible to restrict the
abstract structure of the acting group (Theorem~\ref{last6}). 
In Section~\ref{5} we provide a characterisation of the 
Cartesian decompositions in two of the classes of the 6-Class Theorem
(see Theorems~\ref{c1} and~\ref{diagonal}). 
The investigation of the remaining 4 types can
be found in~\cite{design} and~\cite{1s}. The main results of this
article  are proved in Section~\ref{final}. 

It is also meaningful to study those embeddings of innately transitive
permutation groups $G$ into wreath products $\sym\Gamma\wr\sy\ell$ in
product action in which $G$ induces an intransitive subgroup of
$\sy\ell$. This will be carried out in a separate paper.

Most of our results depend on the correctness of the finite simple
group classification. For instance, a lot of information on the
factorisations of simple and characteristically simple groups that depend
on this classification are used
throughout the paper. Theorem~\ref{disstrips}, which has
central importance in the proof of our 6-Class Theorem and other
theorems in Section~\ref{4}, is a
consequence of~\cite[Lemma~2.2]{bad:quasi}, whose proof uses
the fact, verified by the finite simple group
classification, that each automorphism of a non-abelian, finite simple group
has a non-trivial 
fixed-point. The Schreier Conjecture about the outer automorphism
group of a finite simple group is used in the proof of
Lemma~\ref{minn}.

In this 
paper we use the following notation. Permutations act on the right: if $\pi$ is a permutation and
$\omega$ is a point then the image of $\omega$ under $\pi$ is denoted
$\omega\pi$.   If $G$ is a group acting on a
set $\Omega$ then $G^\Omega$ denotes the subgroup of $\sym\Omega$
induced by $G$.

\section{Further examples of inclusions}\label{exsec}

The examples of inclusions of innately transitive groups into wreath
products in product action given in the introduction are natural. In
this section we give some further examples, whose existence, we believe, 
is really surprising. They indicate the richness of the theory
described in this paper.

Let $\Gamma$ be a finite set, $L\leq\sym\Gamma$, $\ell\geq 2$ an
integer, and $H\leq\sy\ell$. The {\em
wreath product} $L\wr H$ is the semidirect product
$L^\ell\rtimes H$, where, for
$(x_1,\ldots,x_\ell)\in L^\ell$ and $\sigma\in\sy\ell$,
$(x_1,\ldots,x_\ell)^{\sigma^{-1}}=(x_{1{\sigma}},\ldots,x_{\ell{\sigma}})$. The product
action of $L\wr H$ is the action of $L\wr H$ on $\Gamma^\ell$ defined
by
$$
(\gamma_1,\ldots,\gamma_\ell){(x_1,\ldots,x_\ell)}=\left(\gamma_1{x_1},\ldots,\gamma_\ell{x_\ell}\right)\quad\mbox{and}\quad (\gamma_1,\ldots,\gamma_\ell){\sigma^{-1}}=(\gamma_{1\sigma},\ldots,\gamma_{\ell\sigma})
$$
for all $(\gamma_1,\ldots,\gamma_\ell)\in\Gamma^\ell$,
$x_1,\ldots,x_\ell\in L$, and $\sigma\in H$.
The important properties of wreath products can be found in most
textbooks on permutation group theory, see for instance~\cite{dm}.

\begin{example}\label{ex2}
Let $T$ be a finite simple group and let $A,\ B$ be proper subgroups of
$T$ such that $T=AB$. Set $K_1=A\times B$ and $K_2=B\times A$, and let $\Gamma_1$ and $\Gamma_2$ denote the right coset
spaces $[T\times T:K_1]$ and $[T\times T:K_2]$,
respectively. Note that $K_1K_2=T\times T$. Define $\pi\in\sym(\Gamma_1\times\Gamma_2)$ 
as 
$$
(K_1(t_1,t_2),K_2(t_3,t_4))\pi=(K_1(t_4,t_3),K_2(t_2,t_1))$$
for
all $t_1,\ t_2,\ t_3,\ t_4\in T$.
It is routine to check that $\pi$ is an involution and that $\pi$
normalises $\sym\Gamma_1\times\sym\Gamma_2$, interchanging
$\Alt\Gamma_1$ and $\Alt\Gamma_2$. 
Let $\varrho_1$ and $\varrho_2$ be the permutation representations of
$T\times T$ on $\Gamma_1$ and $\Gamma_2$, respectively, 
induced by right multiplication. 
Set
$$
M=\{(\varrho_1(t),\varrho_2(t))\ |\ t\in T\times T\}.
$$
Clearly, $M$ is a subgroup of $\sym\Gamma_1\times\sym\Gamma_2$, and
the stabiliser in $M$ of the point $(K_1,K_2)$ is
$M_0=\{(\varrho_1(t),\varrho_2(t))\ |\ t\in K_1\cap K_2\}$. Since $K_1\cap
K_2=(A\cap B)\times(A\cap B)$, it follows that $|M:M_0|=|T\times
T:K_1\cap K_2|=|T|^2/|A\cap B|^2$. Further, since $T=AB$, we have, for
$i=1,\ 2$, that
$|\Gamma_i|=|T\times T:K_i|=|T|^2/(|A|\cdot|B|)=|T|/|A\cap B|$. Thus
$|M:M_0|=|\Gamma_1\times\Gamma_2|$, and so $M$ is transitive on
$\Gamma_1\times\Gamma_2$.
Let $\alpha=(K_1(t_1,t_2),K_2(t_3,t_4))\in\Gamma_1\times\Gamma_2$ and
$y=(\varrho_1(t),\varrho_2(t))\in M$ with $t=(x_1,x_2)\in T\times T$;
set $y'=(\varrho_1(x_2,x_1),\varrho_2(x_2,x_1))$. Then
\begin{multline*}
\alpha
y\pi=\left(K_1(t_1x_1,t_2x_2)K_2(t_3x_1,t_4x_2)\right)\pi=\left(K_1(t_4x_2,t_3x_1),K_2(t_2x_2,t_1x_1)\right)\\=(K_1(t_1,t_2),K_2(t_3,t_4))\pi
y'=\alpha{\pi
y'}.
\end{multline*}
Since, this is true for all $\alpha\in\Gamma_1\times\Gamma_2$, we have
$y\pi=\pi y'$, and so $y^\pi=y'$ for all $y\in M$. As $y'\in M$, $\pi$ normalises $M$.
This simple argument also shows that
\begin{equation}\label{Mnorm}
\left(\varrho_1(x_1,x_2),\varrho_2(x_1,x_2)\right)^\pi=\left(\varrho_1(x_2,x_1),\varrho_2(x_2,x_1)\right)\quad\mbox{for
all}\quad x_1,\ x_2\in T.
\end{equation}
Set
$$
T_1=\{(\varrho_1(t,1),\varrho_2(t,1))\ |\ t\in T\}\quad\mbox{and}\quad
T_2=\{(\varrho_1(1,t),\varrho_2(1,t))\ |\ t\in T\}.
$$
Then $M$ can be written as the internal direct product $M=T_1\times
T_2$. It also follows from~\eqref{Mnorm} that $\pi$
interchanges $T_1$ and $T_2$. Let
$W=(\sym\Gamma_1\times\sym\Gamma_2)\rtimes\left<\pi\right>$ and
$G=M\rtimes\left<\pi\right>$. Then $G$ is a quasiprimitive
group on $\Gamma_1\times\Gamma_2$ with unique minimal normal subgroup $M$, and $W$ is also quasiprimitive with minimal
normal subgroup $\Alt\Gamma_1\times\Alt\Gamma_2$. Moreover, as $|\Gamma_1|=|\Gamma_2|$,  $W$ is permutationally isomorphic to the
wreath product $\sym\Gamma_1\wr\sy 2$ in product action. Therefore the inclusion $G\leq W$ is as in
case~(c) in Theorem~\ref{main}.
\end{example}

\begin{example}\label{ex3}
Let $T$ be a finite simple group and let $A$, $B$, $C$ be subgroups of
$T$ such that
$T=A(B\cap C)=B(A\cap C)=C(A\cap B)$.
Note that in this case we say that $\{A,B,C\}$ is a strong multiple
factorisation of the finite simple group $T$, and the possibilities
for $T$, $A$, $B$, and $C$ can be found
in~\cite[Table~V]{bad:fact}. Set $K_1=A\times B\times C$, $K_2=B\times
C\times A$, and $K_3=C\times A\times B$; it is a routine calculation
to check that 
\begin{equation}\label{smf}
K_1(K_2\cap K_3)=K_2(K_1\cap K_3)=K_3(K_1\cap
K_2)=T^3.
\end{equation}
Let $\Gamma_1$, $\Gamma_2$, and $\Gamma_3$ denote the right coset
spaces $[T^3:K_1]$, $[T^3:K_2]$, and $[T^3:K_3]$,
respectively, and define $\pi\in\sym(\Gamma_1\times\Gamma_2\times \Gamma_3)$ 
as 
$$
(K_1(t_1,t_2,t_3),K_2(t_4,t_5,t_6),K_3(t_7,t_8,t_9))\pi=(K_1(t_8,t_9,t_7),K_2(t_2,t_3,t_1),K_3(t_5,t_6,t_4))
$$
for
all $t_1,\ldots,t_9\in T$.
It is routine to check that $\pi$ is a permutation of order 3 and that $\pi$
normalises $\sym\Gamma_1\times\sym\Gamma_2\times\sym\Gamma_3$, acting
transitively on $\{\Alt\Gamma_1,\Alt\Gamma_2,\Alt\Gamma_3\}$ by conjugation.
Let $\varrho_1$, $\varrho_2$, and $\varrho_3$ be the permutation representations of
$T^3$ on $\Gamma_1$, $\Gamma_2$, and $\Gamma_3$, respectively, 
induced by right multiplication. 
Set
$$
M=\{(\varrho_1(t),\varrho_2(t),\varrho_3(t))\ |\ t\in T^3\}.
$$
Clearly, $M$ is a subgroup of $\sym\Gamma_1\times\sym\Gamma_2\times
\sym\Gamma_3$. Using~\eqref{smf}, an argument similar to the one in Example~\ref{ex2}
shows that
$M$ is  transitive on $\Gamma_1\times \Gamma_2\times\Gamma_3$. It is
also easy to prove that $\pi$ normalises $M$. 
Set
\begin{eqnarray*}
T_1&=&\{(\varrho_1(t,1,1),\varrho_2(t,1,1),\varrho_3(t,1,1))\ |\ t\in T\},\\ 
T_2&=&\{(\varrho_1(1,t,1),\varrho_2(1,t,1),\varrho_3(1,t,1))\ |\ t\in
T\},\mbox{ and} \\
T_3&=&\{(\varrho_1(1,1,t),\varrho_2(1,1,t),\varrho_3(1,1,t))\ |\ t\in T\}.
\end{eqnarray*}
Then $M$ can be written as the internal direct product $M=T_1\times
T_2\times T_3$, and $\pi$
acts transitively by conjugation on $T_1,\ T_2,\ T_3$. Let 
$W=(\sym\Gamma_1\times\sym\Gamma_2\times\sym\Gamma_3)\rtimes\left<\pi\right>$ and
$G=M\rtimes\left<\pi\right>$. Then $G$ is a quasiprimitive
group with minimal normal subgroup $M$, and $W$ is
quasiprimitive with minimal normal subgroup
$\Alt\Gamma_1\times\Alt\Gamma_2\times\Alt\Gamma_3$. 
Moreover, as $|\Gamma_1|=|\Gamma_2|=|\Gamma_3|$,  
$W$ is 
permutationally isomorphic to  the 
wreath product $\sym\Gamma_1\wr\alt 3$ in product action. Therefore
$G$ can be viewed as a subgroup of $\sym\Gamma_1\wr\sy 3$, and this
inclusion is as in
Theorem~\ref{main}(c).
\end{example}

\begin{example}\label{ex1s}
Let $T$ be a finite simple group and let $A,\ B$ be isomorphic subgroups of $T$
such that $T=AB$. The possibilities for $T$, $A$, and $B$ can be
found in~\cite[Table~2]{recog}. We obtain
from~\cite[Lemma~5.2(ii)]{recog} that there exists a $\vartheta\in\aut T$
that interchanges $A$ and $B$. Define the subgroups $K_1$ and $K_2$ of
$T^4$ by
$$
K_1=\{(t_1,t_2,t_3,t_4)\ |\ t_1\in A,\ t_2\in B,\
t_3=t_4\}\mbox{ and }
K_2=\{(t_1,t_2,t_3,t_4)\ |\ t_1=t_2,\
t_3\in A,\ t_4\in B\}.
$$
One can check using~\cite[Lemma~2.1]{charfact} that $(A\times
B)\{(t,t)\ |\ t\in T\}=T\times T$, and hence $K_1K_2=T^4$. 
Let $\Gamma_1$ and $\Gamma_2$ denote the right coset spaces
$[T^4:K_1]$ and $[T^4:K_2]$, respectively. Define 
$\pi_1\in\sym\Gamma_1\times\sym\Gamma_2$ as follows:
$$
(K_1(t_1,t_2,t_3,t_4),K_2(t_5,t_6,t_7,t_8))\pi_1=\left(K_1\left(\vartheta(t_2),\vartheta(t_1),\vartheta(t_4),\vartheta(t_3)\right),K_2\left(\vartheta(t_6),\vartheta(t_5),\vartheta(t_8),\vartheta(t_7)\right)\right)
$$
for all $t_1,\ldots,t_8\in T$.
Let $\pi_2\in\sym(\Gamma_1\times\Gamma_2)$ be the permutation
defined by
$$
(K_1(t_1,t_2,t_3,t_4),K_2(t_5,t_6,t_7,t_8))\pi_2=(K_1(t_7,t_8,t_5,t_6),K_2(t_3,t_4,t_1,t_2)).
$$
Set $\Pi=\left<\pi_1,\pi_2\right>$. It is routine to check
that $\pi_2$ is an involution that normalises
$\sym\Gamma_1\times\sym\Gamma_2$ swapping the subgroups
$\Alt\Gamma_1$ and $\Alt\Gamma_2$. 
Let $\varrho_1$ and $\varrho_2$ be the permutation representations of
$T^4$ on $\Gamma_1$ and $\Gamma_2$, respectively, 
induced by right multiplication. 
Set
$$
M=\left\{(\varrho_1(t),\varrho_2(t))\ |\ t\in T^4\right\}.
$$
Then $M$ is a subgroup of
$\sym\Gamma_1\times\sym\Gamma_2$. Using the fact that $K_1K_2=T^4$,  it is not
hard to
verify that
$M$ is  transitive on $\Gamma_1\times \Gamma_2$, and it also follows
that $M$ is normalised by $\Pi$. 
Set
\begin{eqnarray*}
T_1&=&\{(\varrho_1(t,1,1,1),\varrho_2(t,1,1,1))\ |\ t\in T\},\\
T_2&=&\{(\varrho_1(1,t,1,1),\varrho_2(1,t,1,1))\ |\ t\in T\},\\
T_3&=&\{(\varrho_1(1,1,t,1),\varrho_2(1,1,t,1))\ |\ t\in T\},\mbox{ and}\\
T_4&=&\{(\varrho_1(1,1,1,t),\varrho_2(1,1,1,t))\ |\ t\in T\}.
\end{eqnarray*}
Then $M$ can be written as the internal direct product $M=T_1\times
T_2\times T_3\times T_4$, and $\Pi$
acts transitively on $T_1,\ T_2,\ T_3,\ T_4$. Let
$W=(\sym\Gamma_1\times\sym\Gamma_2)\rtimes\left<\pi_2\right>$ and
$G=M\Pi$. Then $G$ is a quasiprimitive
group with unique minimal normal subgroup $M$.
Moreover, as $|\Gamma_1|=|\Gamma_2|$,  $W$ quasiprimitive and is
permutationally isomorphic to 
the wreath product $\sym\Gamma_1\wr\sy 2$ in product action. Therefore the inclusion $G\leq W$ is as in
case~(c) in Theorem~\ref{main}.
\end{example}

Note that in all three of the above examples, the crucial fact that made $M$ 
a transitive group on the Cartesian product was that the subgroups $K_i$ formed
a special factorisation of the characteristically simple group $T^2$
in Example~\ref{ex2},
$T^3$ in Example~\ref{ex3}, and $T^4$ in Example~\ref{ex1s}. Crucial
properties of these subgroups $K_i$ are encapsulated in the definition of Cartesian systems in
Section~\ref{cssect}.

\section{Cartesian decompositions and Cartesian systems}\label{cssect}

In a previous paper~\cite{recog} we 
studied the general problem of describing the innately transitive  subgroups of wreath products in product
action. The corresponding problem for primitive groups was solved
by~\cite{prae:inc} and~\cite{blowups}, but for quasiprimitive groups
it was left open in~\cite{bad:quasi}. Related
problems were also addressed
in~\cite{baum}. 
If an innately transitive group $G$
with plinth $M$ is a subgroup of such a wreath product $W$, then
the underlying set can be viewed as a Cartesian product of
smaller sets. The actions of the groups $W$ and $G$ preserve this
Cartesian product. We made these ideas more precise by introducing the
concept of a Cartesian decomposition of a set.

A
{\em Cartesian decomposition} of a set $\Omega$ is a set
$\{\Gamma_1,\ldots,\Gamma_\ell\}$ of partitions of $\Omega$ such that 
$$
|\gamma_1\cap\cdots\cap\gamma_\ell|=1\quad\mbox{for
all}\quad\gamma_1\in\Gamma_1,\ldots,\gamma_\ell\in\Gamma_\ell.
$$
This property implies that the map
$\omega\mapsto(\gamma_1,\ldots,\gamma_\ell)$, where for
$i=1,\ldots,\ell$ the block $\gamma_i\in\Gamma_i$ is chosen so that
$\omega\in\gamma_i$, is a well defined bijection between $\Omega$ and
$\Gamma_1\times\cdots\times \Gamma_\ell$. Thus the set $\Omega$ can
naturally be
identified with the Cartesian product
$\Gamma_1\times\cdots\times\Gamma_\ell$. This definition was first suggested
by~\cite{kov:decomp}.

If $G$ is a permutation group acting on $\Omega$, then a Cartesian
decomposition $\E$ of $\Omega$ is said to be $G$-invariant if the partitions in
$\E$ are permuted by $G$.
For a permutation group $G\leq\sym\Omega$, the symbol $\cd {}G$
denotes the set of $G$-invariant Cartesian decompositions of
$\Omega$. If $\E\in\cd {}G$ and $G$ acts on $\E$
transitively, then $\E$ is said to be a {\em transitive} $G$-invariant
Cartesian decomposition. The set of transitive $G$-invariant
Cartesian decompositions of $\Omega$ is denoted $\cd{\rm
  tr}G$. 

If $L\leq\sym\Gamma$ and $H\leq\sy\ell$ with some $\ell\geq 2$, 
then the wreath product $W=L\wr H$ at the beginning of Section~\ref{exsec} is
considered as a permutation group acting in product action on the set $\Gamma^\ell$. There is a
natural Cartesian decomposition of $\Gamma^\ell$ preserved by
$W$; namely, take $\E=\{\Gamma_1,\ldots,\Gamma_\ell\}$ where
$$
\Gamma_i=\{\{(\gamma_1,\ldots,\gamma_\ell)\ |\ \gamma_i=\gamma\}\ |\
\gamma\in\Gamma\}\quad\mbox{for}\quad i=1,\ldots,\ell. 
$$
The reader can easily check that $\E$ is a
$W$-invariant Cartesian decomposition of $\Gamma^\ell$, and that the
$W$-actions on the set $\E$ and on the set of natural coordinate subgroups of
the base group $L^\ell$ of $W$ are equivalent. We study innately
transitive subgroups of wreath products in product action, such as
$W$, via the natural Cartesian decomposition of the underlying set
corresponding to the product action of $W$.

Suppose that $G$ is an innately transitive
subgroup of $\sym\Omega$ with plinth $M$, and that $\E$ is a
$G$-invariant Cartesian decomposition of $\Omega$. 
We proved
in~\cite[Proposition~2.1]{recog} that each $\Gamma_i\in\E$ is an
$M$-invariant partition of $\Omega$. 
Choose an element $\omega$ of $\Omega$ and let 
$\gamma_1\in\Gamma_1,\ldots,\gamma_\ell\in\Gamma_\ell$ be such that
$\{\omega\}=\gamma_1\cap\cdots\cap\gamma_\ell$; set $K_i=M_{\gamma_i}$. 
Then
\cite[Lemmas~2.2 and 2.3]{recog} imply that the set $\K_\omega(\E)=\{K_1,\ldots,K_\ell\}$ is invariant
under conjugation by $G_\omega$, and moreover
\begin{eqnarray}\label{csdef1}
\bigcap_{i=1}^\ell K_i&=&M_\omega\quad\mbox{and}\\
\label{csdef2}K_i\left(\bigcap_{j\neq
i}K_j\right)&=&M\quad \mbox{for all}\quad i\in\{1,\ldots,\ell\}.
\end{eqnarray}

For an arbitrary transitive permutation group $M$ on $\Omega$, and a
point $\omega\in\Omega$,
a set $\K=\{K_1,\ldots,K_\ell\}$ of proper subgroups
of $M$ is said to be 
a {\em Cartesian system of subgroups with respect to $\omega$} 
for $M$, if~\eqref{csdef1} and \eqref{csdef2}~hold. If
$M$ is an abstract group then a set $\{K_1,\ldots,K_\ell\}$ of
proper subgroups satisfying~\eqref{csdef2} is said to be a {\em Cartesian
system}.

The reader can easily see that the subgroups $K_i$ in
Examples~\ref{ex2}, \ref{ex3}, and \ref{ex1s} form, in each example, 
a Cartesian system 
for the characteristically simple group $M$. These examples show the
importance of Cartesian systems for constructing inclusions of an
innately transitive group $G$ into wreath products in product action. 
In fact, Cartesian
systems provide an effective way of identifying the set of $G$-invariant Cartesian decompositions from
information internal to $G$. 

\begin{theorem}[Theorem~1.4 and Lemma~2.3~\cite{recog}]\label{bij}
Let $G\leq \sym\Omega$ be an innately transitive permutation group with
plinth $M$. For a fixed $\omega\in\Omega$ 
the correspondence $\E\mapsto\K_\omega(\E)$ is a bijection between the
set of $G$-invariant Cartesian decompositions of $\Omega$
and the set of $G_\omega$-invariant
Cartesian systems of subgroups for $M$ with respect to
$\omega$. Moreover the $G_\omega$-actions on $\E$ and on $\K_\omega(\E)$
are equivalent.
\end{theorem}

Suppose that $G\leq\sym\Omega$ is an innately transitive group with
plinth $M$, and let $\omega\in\Omega$. Let $\K$ be a
$G_\omega$-invariant  Cartesian system of subgroups in $M$ with
respect to $\omega$. Then the
previous theorem implies that $\K=\K_\omega(\E)$ for some 
$G$-invariant Cartesian decomposition $\E$ of $\Omega$. Indeed, it is easy
to see that $\E$ consists of the $M$-invariant partitions 
$\{(\omega^K)^m\ |\ m\in M\}
$
where $K$ runs through the elements of $\K$. 
This Cartesian decomposition
is denoted $\E(\K)$.

Using this
theory we were able to describe in~\cite{recog} those innately transitive subgroups 
of  wreath products that have a simple plinth. This led to
a classification of transitive simple and almost simple subgroups of
wreath products in product action
(see~\cite[Theorem~1.1]{recog}). Our aim in this paper is to extend that theory
to achieve a more complete characterisation of innately transitive
subgroups of wreath products in product action for innately transitive
groups with a non-abelian plinth. As noted earlier, the case of
abelian plinth was settled in~\cite{prae:inc}.

\section{On characteristically simple groups}\label{charsec}

Property~\eqref{csdef2} is of crucial importance for investigating Cartesian systems in
characteristically simple groups. 
Hence it is very important for our research
to study factorisations of such groups. In
addition to the results in~\cite{charfact}, we use the ones listed in this section.

Suppose that $T$ is a finite, non-abelian simple group and $\L$ is a
set of proper subgroups in $T$ such that $|\L|\geq 3$ and 
$A(B\cap C)=T$ whenever $A$,
$B$, and  $C$ are pairwise distinct elements of $\L$. Then the set $\L$ is said to be a {\em strong multiple
factorisation} of $T$. Strong multiple factorisations of finite simple
groups were classified in~\cite[Table~V]{bad:fact}. It was proved that
not all finite
simple groups admit a strong multiple factorisation, and each such
factorisation contains exactly three, pairwise non-isomorphic subgroups.

Suppose that $M=T_1\times\cdots\times T_k$ where the $T_i$ are
non-abelian, 
finite
simple groups. For
$I\subseteq\{T_1,\ldots,T_k\}$ the function $\sigma_I:M\rightarrow\prod_{T_i\in
I}T_i$ is the natural projection map. We also write
$\sigma_i$ for $\sigma_{\{T_i\}}$. 
A subgroup $X$ of $M$ is said to be a {\em strip} if for
each $i=1,\ldots,k$ either $\sigma_i(X)=1$ or $\sigma_i(X)\cong X$. The set
of $T_i$ such that $\sigma_i(X)\neq 1$ is called the {\em support} of
$X$ and is denoted $\supp X$. If $T_m\in\supp X$ then we also say that
$X$ {\em covers} $T_m$. Two strips $X_1$ and $X_2$ are {\em disjoint}
if $\supp X_1\cap\supp X_2=\emptyset$. A strip $X$ is said to be {\em full} if
$\sigma_i(X)=T_i$ for all $T_i\in\supp X$,
and $X$ is called {\em non-trivial} if
$|\supp X|\geq 2$. A subgroup $K$ of $M$ is said to be {\em subdirect} if
$\sigma_i(K)=T_i$ for all $i$.

We recall a well-known lemma on finite simple groups which can be
found in~\cite{scott}. 

\begin{lemma}\label{scott}
Let $M$ be a direct product of finitely many non-abelian, finite simple groups and $H$ a
subdirect subgroup of $M$. Then $H$ is the direct product of pairwise
disjoint full strips of $M$.
\end{lemma}

The following result gives a generalisation of Scott's
Lemma.  Let $M=T_1\times\cdots\times T_k$ be a characteristically
simple group where $T_1,\ldots,T_k$ are the simple normal subgroups of
$M$. If $K$ is a subgroup of $M$ and $X$ is a non-trivial strip in
$M$ such that $K=X\times\sigma_{\{T_1,\ldots,T_k\}\setminus\supp X}(K)$ then we say that
$X$ is {\em involved} in $K$.

\begin{lemma}\label{scottgen}
Let $M$ be a direct product $T_1\times\cdots\times T_k$ where the
$T_i$ are pairwise isomorphic finite simple groups, and let $K$ be a subgroup of
$M$ such that
$\sigma_m(K)=T_m$ for some $m\leq k$. Then there is a unique full strip $X$ of
$M$ covering $T_m$ such that $X$ is involved in $K$.  
\end{lemma}
\begin{proof}
Let $X$, $Y$ be normal subgroups of $K$ which are minimal by inclusion
subject to satisfying $\sigma_m(X)=\sigma_m(Y)=T_m$; such  
subgroups exist as $\sigma_m(K)=T_m$. 
Then
\[\sigma_m(X\cap Y)\geq\sigma_m([X,Y])=\left[\sigma_m(X),\sigma_m(Y)\right]=[T_m,T_m]=T_m\]
and by minimality we have $X=X\cap Y=Y$. Thus $K$ has a
unique normal subgroup, $X$ say, which is minimal by inclusion subject to
satisfying $\sigma_m(X)=T_m$.

We show that $X\cong T_1$.
For $t\not=m$ consider $\sigma_t(X)$; suppose that this is
non-trivial whence $X\cap \ker \sigma_t$ is a proper
subgroup of $X$.
As $X\cap\ker \sigma_t$ is normal in $X$ we deduce that
$\sigma_m(X\cap \ker \sigma_t)$ is a normal subgroup
of $T_m$ and is not equal to $T_m$ by the 
minimality of $X$. Thus $\sigma_m(X\cap \ker\sigma_t)$ is trivial,
and so $X\cap \ker \sigma_t\leq X\cap \ker \sigma_m$. The
latter
subgroup has index $|T_m|$ in $X$, while the former has index
dividing $|T_t|=|T_m|$. This forces
$X\cap\ker\sigma_m=X\cap\ker\sigma_t$ and
$\sigma_t(X)=T_t$. So for each $t=1,\ldots,k$ we have either
$\sigma_t(X)=1$, or $\sigma_t(X)=T_t$ and $X\cap\ker
\sigma_t=X\cap\ker\sigma_m$. Hence $\sigma_t(X\cap
\ker\sigma_m)=1$ for all $t$, whence
$X\cap\ker\sigma_m$ is trivial and $X\cong T_m\cong T_1$.

Thus $X$ is a full strip which covers $T_m$.
As $X$ is
self-normalising in $\prod_{T_t\in\supp X}T_t$ (see~\cite[Lemma~4.6]{bp}),
we see that $\sigma_{\supp X}(K)=X$ whence
$K=X\times\sigma_{\{T_1,\ldots,T_k\}\setminus\supp X}(K)$. Hence $X$
is involved in $K$.
\end{proof}

The next lemma implies that the factorisation of a characteristically
simple group as a product of a full strip and a proper subgroup is
possible only under restricted conditions. This result is crucial for
the proof of Theorem~\ref{stripth}.

\begin{lemma}\label{fsf}
Suppose that $M=T_1\times\cdots\times T_k$ is a non-abelian 
characteristically simple group, $\alpha_i:T_1\rightarrow T_i$ is  an isomorphism for $i=2,\ldots,k$,
and $K$ is a subgroup of $M$ such
that $\sigma_i(K)\neq T_i$ for $i=1,\ldots,k$. If 
$$
\{(t,\alpha_2(t),\ldots,\alpha_k(t))\ |\ t\in T_1\}K=M,
$$
then
$k\leq 3$. 
Moreover, if $k=3$ then
$\{\sigma_1(K),\alpha_2^{-1}(\sigma_2(K)),\alpha_3^{-1}(\sigma_3(K))\}$
is a strong multiple factorisation of $T_1$. 
\end{lemma}
\begin{proof}
For each $t\in T_1$
there is some $x\in T_1$, and $(a_1,\ldots,a_k)\in K$ such that
$$
(t,1,\ldots,1)=(x,\alpha_2(x),\ldots,\alpha_k(x))(a_1,\ldots,a_k).
$$ 
Then $t=xa_1$ and $x=\alpha_i^{-1}(a_i^{-1})$ for all $i=2,\ldots,k$, and so
$x\in\alpha_2^{-1}(\sigma_2(K))\cap\cdots\cap\alpha_k^{-1}(\sigma_k(K))$. So
$$
T_1=\sigma_1(K)(\alpha_2^{-1}(\sigma_2(K))\cap\cdots\cap\alpha_k^{-1}(\sigma_k(K))),
$$
and we obtain similarly that
\begin{multline*}
T_1=\alpha_2^{-1}(\sigma_2(K))(\sigma_1(K)\cap\alpha_3^{-1}(\sigma_3(K))\cap\cdots\cap
\alpha_k^{-1}(\sigma_k(K)))\\=\cdots=\alpha_k^{-1}(\sigma_k(K))(\sigma_1(K)\cap\alpha_2^{-1}(\sigma_2(K))\cap\cdots\cap\alpha_{k-1}^{-1}(\sigma_{k-1}(K))).
\end{multline*}
Since $\sigma_i(K)$ is a proper subgroup of $T_i$ for all $i$, we
obtain that $$
\{\sigma_1(K),\alpha_2^{-1}(\sigma_2(K)),\ldots,\alpha_k^{-1}(\sigma_k(K))\}$$
is a strong multiple factorisation of $T_1$ whenever $k\geq 3$. 
Then the results of
Baddeley and Praeger~\cite{bad:fact} imply that $k=3$. 
\end{proof}

The next lemma which can be found as~\cite[Lemma~2.2]{bad:quasi} states that finite
characteristically simple groups cannot be written as a product of two
subgroups each of which is a direct product of non-trivial
strips. The proof of this result uses the fact that each automorphism
of a non-abelian, finite simple group has a non-trivial fixed-point,
and hence it depends on the finite simple group classification.

\begin{lemma}\label{2.1}
Suppose that $M=T_1\times \cdots\times T_k$ is a direct product of
isomorphic non-abelian, simple groups $T_1,\ldots,T_k$. Suppose that
$A_1,\ldots,A_{m}$ are non-trivial pairwise disjoint strips in $M$,
and so are $B_1,\ldots,B_n$. Then $M\neq (A_1\times\cdots\times
A_m)(B_1\times\cdots\times B_n)$. 
\end{lemma}

\section{Strips in a Cartesian system}\label{stripsec}

A non-trivial strip $X$ is said to be involved in a Cartesian system
$\K$ for a non-abelian, characteristically simple group
if $X$ is involved in an element of $\K$, as defined
before Lemma~\ref{scottgen}. 
Note that in this case
Lemma~\ref{2.1} and~\eqref{csdef2} 
imply that $X$ is involved in a unique element of $\K$. The purpose of
this section is to prove that two distinct strips involved in a Cartesian system
must be disjoint. This result plays a key r\^ole in the proof of our
6-Class Theorem (see in particular the proof of Theorem~\ref{stripth}(d)).

\begin{theorem}\label{disstrips}
Let $M=T_1\times\cdots\times T_k$ be a finite, non-abelian, characteristically simple
group with simple normal subgroups $T_1,\ldots,T_k$, and let $G_0$ be a subgroup of $\aut M$ such that the natural
action of $G_0$ on $T_1,\ldots,T_k$ is transitive.
Suppose, in addition, that $\K=\{K_1,\ldots,K_\ell\}$ is a
$G_0$-invariant Cartesian system of subgroups in the abstract group
$M$. 
If $X_1$ and $X_2$ are distinct, non-trivial strips involved in $\K$, then $X_1$ and $X_2$ are disjoint.
\end{theorem}
\begin{proof}
By the definition given above,  if $X_1$, $X_2$ are involved in the same
$K_i$ then
they are disjoint as strips. Thus we may assume that $X_1$ is involved
in $K_{j_1}$ and $X_2$ is involved in $K_{j_2}$ where $j_1\neq
j_2$. First we prove that $|\supp X_1\cap\supp X_2|\leq 1$. Suppose to
the contrary that $T_{i_1},\
T_{i_2}\in\supp X_1\cap \supp X_2$ with $i_1\neq i_2$. Then
$\sigma_{\{T_{i_1},T_{i_2}\}}(K_{j_1})$ and
$\sigma_{\{T_{i_1},T_{i_2}\}}(K_{j_2})$ are non-trivial strips, and so Lemma~\ref{2.1} implies that
$\sigma_{\{T_{i_1},T_{i_2}\}}(K_{j_1})\sigma_{\{T_{i_1},T_{i_2}\}}(K_{j_2})\neq
T_{i_1}\times T_{i_2}$. However $K_{j_1}K_{j_2}=M$ by~\eqref{csdef2}, 
which is a
contradiction. Thus $|\supp X_1\cap\supp X_2|\leq 1$. 

Assume that $\supp{X_1}\cap \supp{X_2}=\{T_t\}$ for some $t\leq k$.
Choose $g\in G_0$ such that $T_t^g\in \supp{X_2}\setminus
\supp{X_1}$; such an element $g$ exists since $G_0$ is transitive on
$\T=\{T_1,\ldots,T_k\}$ and $\supp{X_2}\setminus \supp{X_1}$ is non-empty.
Now $G_0$ acts by conjugation on the set of strips involved in $\K$,
and so both $X_1^g$ and $X_2^g$ are  strips involved in $\K$. As
$T_t^g$ is in both $\supp{X_1^g}$ and $\supp{X_2^g}$, but is not in
$\supp{X_1}$ we deduce that there exists a non-trivial  strip $X_3$ in $\K$
distinct from $X_1$, $X_2$ such that $\supp{X_3}\cap\left(\supp{X_2}\setminus
\supp{X_1}\right)\not=\emptyset$
(namely, we can take $X_3$ to be one of $X_1^g$ or $X_2^g$ as at least one
of these is distinct from $X_1$ and $X_2$).
Proceeding in this way we construct a sequence
$X_1,X_2,\ldots$ of distinct, non-trivial strips in $\K$ such
that $\supp{X_{d+1}}\cap\left(\supp{X_d}\setminus
\supp{X_{d-1}}\right)\not=\emptyset$
for each $d\geq 2$.
Let $X_a$ be the first member of the sequence with $a\geq3$
and
\[\supp{X_a}\cap\bigl(\supp{X_1}\cup\cdots\cup
\supp{X_{a-2}}\bigr)\not=\emptyset. \]
By removing some initial segment of the sequence if
necessary
we may assume that the intersection $\supp{X_a}\cap \supp{X_1}$ is
non-empty,
while $\supp{X_a}\cap \supp{X_d}=\emptyset$ if $2\leq d \leq a-2$ for some $a\geq 3$.

By relabeling the $K_s$ we may assume that $X_1$ is involved in $K_1$.
Let $1=i_1<i_2<\cdots<i_d<a$ be such that among the $X_i$ 
the strips $X_{i_1},\ldots,X_{i_d}$ are precisely the ones
that are involved in $K_1$. 
Note that $X_a$ is not involved in $K_1$ since $\supp{X_a}$ and $\supp{X_1}$ are not
disjoint. Also, $i_{j+1}\geq i_j+2$ for all
$j=1,\ldots,d-1$
since $\supp{X_{i_j}}$ and $\supp{X_{i_j+1}}$ are not disjoint. We
may also relabel the $T_i$ so that
$$
\{T_1\}=\supp{X_a}\cap
\supp{X_1}\quad\mbox{and}\quad\{T_2\}=\supp{X_1}\cap \supp{X_2},
$$
and so that for $j=2,\ldots,d$,
$$\{T_{2j-1}\}=\supp{X_{i_{j}-1}}\cap \supp{X_{i_j}}\quad\mbox{and}\quad
\{T_{2j}\}=\supp{X_{i_j}}\cap \supp{X_{i_{j}+1}}.
$$
It follows from the definition of $a$ that $T_1,\ldots,T_{2d}$ are
pairwise distinct.
Let $\sigma$ be the projection map $M\to
T_1\times\cdots\times T_{2d}$. The labeling of $T_1,\ldots,
T_{2d}$ and the minimality of $a$ have ensured that
the following all hold:
\begin{enumerate}
\item[(i)] $\sigma\left(X_{i_j}\right)$ is a non-trivial strip of $T_{2j-1}\times T_{2j}$ for
$j=1,\ldots,d$;
\item[(ii)] for $j=1,\ldots,d-1$ we have
$$
\sigma\left(X_{i_j+1}\cap
\cdots\cap X_{i_{j+1}-1}\right)=
\left\{\begin{array}{ll}
1 & \mbox{if $i_{j+1}\geq i_j+3$}\\
\mbox{a non-trivial strip covering $T_{2j}$ and $T_{2j+1}$} & \mbox{if
$i_{j+1}=i_j+2$};\end{array}\right.
$$
\item[(iii)] moreover,
$$
\sigma\left(X_{i_d+1}\cap
\cdots\cap X_{a}\right)=
\left\{\begin{array}{ll}
1 &\mbox{if $a\geq i_d+2$}\\ 
\mbox{a non-trivial strip covering $T_{2d}$ and $T_1$}&\mbox{if $a=i_d+1$.}\end{array}\right.
$$
\end{enumerate}
Recalling that $X_{i_1},\ldots,X_{i_d}$ are precisely the
strips from the sequence $X_1,\ldots,X_a$ that appear as
direct factors of $K_1$, we have that
\[\sigma(K_1)=\sigma(X_{i_1}\times\cdots\times X_{i_d}),\]
and that for $\widehat K_1=K_2\cap \cdots\cap K_\ell$,
\[\sigma(\widehat K_1)\leq \sigma\Bigl(
(X_2\cap\cdots\cap X_{i_2-1})\times \cdots\times
(X_{i_d+1}\cap\cdots\cap X_{a})\Bigr).\]
From (i) we deduce that $\sigma(K_1)$ is the direct
product of strips of $T_1\times
T_2,\ldots,
T_{2d-1}\times T_{2d}$, while from (ii) and (iii)
we deduce that $\sigma(
\widehat K_1)$ is contained in the direct
product of non-trivial strips of $T_2\times
T_3,\ldots,T_{2d-2}\times T_{2d-1},T_{2d}\times T_1$.
Thus $\sigma(K_1)$ and $\sigma(\widehat K_1)$ are each contained
in the
direct product of disjoint non-trivial strips of $T_1\times
\cdots\times T_{2d}$.
Lemma~\ref{2.1} implies  that
$\sigma(K_1)\sigma(\widehat K_1)\not=\sigma(M)$,
which contradicts the fact implied by~\eqref{csdef2} that $K_1\widehat K_1=M$. Thus $\supp
X_1\cap\supp X_2=\emptyset$, as required.
\end{proof}

\section{Six classes of Cartesian decompositions}\label{4}

A non-abelian plinth of an innately transitive group 
$G$ has the form $M=T_1\times\cdots\times T_k$ where
the $T_i$ are  finite, non-abelian, simple groups. Let $\E\in\cd{}G$ and let
$\K_\omega(\E)$ be a corresponding Cartesian system
$\{K_1,\ldots,K_\ell\}$ for $M$ with respect to some $\omega\in\Omega$. Then equation~\eqref{csdef2} implies that, for all $i\leq k$ and $j\leq\ell$,
\begin{equation}\label{simpfact}
\sigma_i(K_j)\left(\bigcap_{j'\neq j}\sigma_i(K_{j'})\right)=T_i.
\end{equation}
In particular this means that if $\sigma_i(K_j)$ is a proper subgroup
of $T_i$ then $\sigma_i(K_{j'})\neq\sigma_i(K_j)$ for all
$j'\in\{1,\ldots,\ell\}\setminus\{j\}$. 
It is thus important to understand the following sets of subgroups:
\begin{equation}\label{f}
\mathcal F_i(\E,M,\omega)=\{\sigma_i(K_j)\ |\ j=1,\ldots,\ell,\ \sigma_i(K_j)\neq T_i\}.
\end{equation}
Note that $|\mathcal F_i(\E,M,\omega)|$ is the number of indices $j$ such that
$\sigma_i(K_j)\neq T_i$. 
The set $\mathcal F_i(\E,M,\omega)$ is independent of $i$ up to isomorphism, 
in the sense that if $i_1,\ i_2\in\{1,\ldots,k\}$ and $g\in G_\omega$ 
are such that
$T_{i_1}^g=T_{i_2}$ then $\mathcal
F_{i_1}(\E,M,\omega)^g=\{L^g\ |\ L\in\mathcal F_i(\E,M,\omega)\}=\mathcal F_{i_2}(\E,M,\omega)$. This argument
also shows that the elements of $\mathcal F_{i_1}(\E,M,\omega)$ are actually 
$G_\omega$-conjugate to the elements of $\mathcal F_{i_2}(\E,M,\omega)$.

Recall that, for $G\leq\sym\Omega$, $\cd{\rm tr}G$ denotes the set of
transitive $G$-invariant Cartesian decompositions 
of $\Omega$. The following theorem depends on
the finite simple group classification, since the proof uses results
from~\cite{bad:fact} on full factorisations and multiple
factorisations of finite simple groups.

\begin{theorem}\label{stripth}
Suppose that $G$ is an innately transitive permutation group with
a non-abelian plinth $M=T_1\times\cdots\times T_k$ where $T_1,\ldots,T_k$ are
pairwise isomorphic finite simple groups and $k\geq 1$. Let $\E\in\cd {\rm tr}G$
with a corresponding Cartesian system $\K$ for $M$ with respect to
$\omega\in\Omega$. For $i=1,\ldots,k$, let
$\mathcal F_i=\mathcal F_i(\E,M,\omega)$ be defined as
in~\eqref{f}. Then the following all hold.
\begin{enumerate}
\item[(a)] The number $|\mathcal F_i|$ is independent of $i$ and
$|\mathcal F_i|\leq 3$. Further, if $|\mathcal F_i|=3$, then $\mathcal
F_i$ is a strong multiple factorisation of $T_i$.
\item[(b)] Suppose that there is a non-trivial, full strip involved in
$\K$. Then $k\geq 2$ and $|\mathcal F_i|\in\{0,\ 1\}$. Moreover, if
$|\mathcal F_i|=1$, then the $T_i$ admit a
factorisation as a product of two proper, isomorphic subgroups.
\item[(c)] If $X$ is a non-trivial, full strip involved in $\K$ and
$|\mathcal F_i|=1$ then $|\supp X|=2$. 
\item[(d)] Set
$\mathcal P=\{\supp X\ |\ X\mbox{ is a non-trivial, full strip
involved in $\K$}\}$. If $\mathcal P\neq\emptyset$ then $\mathcal P$
is a $G$-invariant partition of $\{T_1,\ldots,T_k\}$. 
\end{enumerate}
\end{theorem}
\begin{proof}
(a) Suppose that $i_1,\ i_2\in\{1,\ldots,k\}$ and $g\in G_\omega$ such
that $T_{i_1}^g=T_{i_2}$. Then it was observed after equation~\eqref{f}
that $\mathcal F_{i_1}^g=\mathcal F_{i_2}$, and so $|\mathcal F_i|$ is
independent of $i$. The definition of $\mathcal F_i$
and~\eqref{simpfact} imply that if $|\mathcal F_i|\geq 3$ then 
$\mathcal F_i$ is a strong multiple factorisation of
$T_i$. Hence~\cite[Table~V]{bad:fact} shows that $|\mathcal F_i|\leq
3$. 

(b) A non-trivial strip has to cover at least two of the $T_i$, and so
we must have $k\geq 2$. We argue by contradiction and assume that $|\mathcal F_i|\geq 2$.
Suppose without loss of generality that $K_1$ involves a non-trivial full strip $X$
covering $T_1$ and $T_2$. Let $g\in G_\omega$ such that $T_1^g=T_2$. 
Then, by Theorem~\ref{disstrips}, $X^g=X$, and so $K_1^g=K_1$, and also
$\widehat K_{1}^g=\widehat K_{1}$, where, as in the proof of
Theorem~\ref{disstrips}, $\widehat K_1=K_2\cap\cdots\cap K_\ell$.
This also implies that $\sigma_1(\widehat K_1)^g=\sigma_2(\widehat
K_1)$; therefore $\sigma_1(\widehat K_1)$ and $\sigma_2(\widehat K_1)$
are isomorphic. Note that
$$
\sigma_{\{1,2\}}(X)=\{(t,\alpha(t))\ |\ t\in T_1\}
$$
for some isomorphism $\alpha:T_1\rightarrow T_2$. As $|\mathcal
F_1|\geq 2$, we obtain that $\sigma_1(K_{j_1})<T_1$ and
$\sigma_1(K_{j_2})<T_1$ for distinct $j_1,\ j_2\in\{2,\ldots,\ell\}$,
and, as noted above, $\sigma_1(K_{j_1})\neq\sigma_1(K_{j_2})$.  Thus
$\sigma_1(\widehat K_1)<T_1$ and, as $\sigma_1(\widehat
K_1)\cong\sigma_2(\widehat K_1)$, also  $\sigma_2(\widehat K_1)<T_2$. Since
$K_1\widehat K_{1}=M$, \cite[Lemma 2.1]{charfact} implies that
$\sigma_1(\widehat K_{1})$ and $\alpha^{-1}(\sigma_2(\widehat K_{1}))$ form a
full factorisation of $T_1$ with isomorphic subgroups. Based on the
classification of full factorisations of almost simple groups
in~\cite{bad:fact}, such factorisations were 
described in~\cite[Lemma~5.2]{recog}. In particular,
$\sigma_1(\widehat K_{1})$ and $\alpha^{-1}(\sigma_2(\widehat K_{1}))$ are
both maximal subgroups of $T_1$.
On the other hand,
$\sigma_1(\widehat K_{1})\leq
\sigma_{1}(K_{j_1})\cap\sigma_1(K_{j_2})<\sigma_1(K_{j_1}),\
\sigma_1(K_{j_2})$, which is a contradiction. Therefore $|\mathcal
F_i|\leq 1$. 

(c) 
Suppose without loss of generality that $X$ is a non-trivial, full strip involved in
$K_1$ covering
$T_1,\ldots,T_s$. Let $i,\ j$ be distinct elements of $\{1,\ldots,s\}$
and let $g\in G_\omega$ such that $T_i^g=T_j$. Then $X^g$ is a strip involved in $\K$ such that $T_j\in\supp X\cap
\supp X^g$. Therefore by Theorem~\ref{disstrips}, $X=X^g$, and so
$K_1^g=K_1$, and also $\widehat K_{1}^g=\widehat K_{1}$. 
Since $T_i^g=T_j$ it follows that
$\sigma_i(\widehat K_1)^g=\sigma_j(\widehat K_1)$. Now, as
$|\mathcal F_i|=1$, there exists $t$ such that $\sigma_i(K_t)<T_i$,
and, as $\sigma_i(K_1)=T_i$, we have $t\neq 1$. Hence 
$\sigma_i(\widehat K_1)<T_i$. Also,
$\sigma_j(\widehat K_1)=\sigma_i(\widehat K_1)^g<T_i^g=T_j$.
Therefore $\sigma_i(\widehat K_1)<T_i$ holds for all
$i\in\{1,\ldots,s\}$.
Thus
the factorisation $\sigma_{\supp X}(M)=\sigma_{\supp X}(K_1)\sigma_{\supp
X}(\widehat K_{1})$ is as in Lemma~\ref{fsf}, and so $s\leq 3$, and moreover if
$s=3$, then the projections $\sigma_1(\widehat K_1)$,
$\sigma_2(\widehat K_1)$, and $\sigma_3(\widehat K_1)$ are isomorphic
to the subgroups in  a strong multiple factorisation of $T_1$. On the
other hand, the subgroups in such a strong multiple factorisation
are pairwise non-isomorphic (see~\cite[Table~V]{bad:fact}): a contradiction. Hence we obtain $s=2$. 

(d) Let $\mathcal P=\{\supp X\ |\ X\mbox{ is a non-trivial, full strip
involved in $\K$}\}$. By Theorem~\ref{disstrips}, either $\mathcal
P=\emptyset$ or $\mathcal P$ 
is a partition of $\{T_1,\ldots,T_k\}$. Moreover,
if $\supp X\in\mathcal P$ and $g\in G_\omega$, then there exists
$K\in\K$ such that $X$ is involved in $K$. Therefore $K^g\in\K$ and
$X^g$ is involved in $K^g$. Thus $X^g$ is involved in $\K$, and so
$(\supp X)^g=\supp (X^g)\in\mathcal P$. Hence $\mathcal P$ is
$G_\omega$-invariant. Since $M$ acts trivially on $\{T_1,\ldots,T_k\}$
by conjugation and $G=MG_\omega$, we have that $\mathcal P$ is a
$G$-invariant partition of $\{T_1,\ldots,T_k\}$.
\end{proof}

If $G$ is a finite,  innately transitive group with a non-abelian plinth $M$,
then the set $\cd{\rm tr}G$ is further subdivided according to
the structure of the subgroups in the corresponding Cartesian systems
for $M$ as follows. For $\E\in\cd{\rm tr}G$ and $\omega\in\Omega$ the
sets $\mathcal F_i=\mathcal F_i(\E,M,\omega)$ are as defined in~\eqref{f}.
\begin{eqnarray*}
\cd{\rm S}G&=&\{\E\in\cd{\rm tr}G\ |\ \mbox{the elements of $\K_\omega(\E)$
are subdirect subgroups in $M$}\};\\
\cd 1G&=&\{\E\in\cd{\rm tr}G\ |\ |\mathcal F_i|=1\mbox{ and $\K_\omega(\E)$
involves no non-trivial, full strip}\};\\
\cd {\rm 1S}G&=&\{\E\in\cd{\rm tr}G\ |\ |\mathcal F_i|=1\mbox{ and $\K_\omega(\E)$
involves non-trivial, full strips}\};\\
\cd {2\sim}G&=&\{\E\in\cd{\rm tr}G\ |\ |\mathcal F_i|=2\mbox{ and
the $\mathcal F_i$ contain two $G_\omega$-conjugate subgroups}\};\\
\cd{2\not\sim}G&=&\{\E\in\cd{\rm tr}G\ |\ |\mathcal F_i|=2\mbox{ and
the subgroups in $\mathcal F_i$ are not $G_\omega$-conjugate}\};\\
\cd 3G&=&\{\E\in\cd{\rm tr}G\ |\ |\mathcal F_i|=3\}.
\end{eqnarray*}

In order to prove that the above classes of Cartesian decompositions
are well-defined, we need to show that the properties used in these
definitions are independent of the choice of the plinth $M$, the
Cartesian system $\K_\omega(\E)$, and the index $i$. 

\begin{theorem}[6-class Theorem]\label{5class}
If $G$ is a finite, innately transitive permutation group with a non-abelian
plinth $M$, then the classes $\cd 1G$, $\cd SG$, $\cd {1S}G$,
$\cd{2\sim}G$, $\cd {2\not\sim}G$, and $\cd 3G$
are well-defined subsets of $\cd{\rm tr}G$, and they form a partition of $\cd{\rm tr}G$. Moreover, if $M$ is simple
then $\cd{\rm tr}G=\cd {2\sim}G$. 
\end{theorem}
\begin{proof}
First we prove that the above classes are well-defined. 
Suppose that $\E\in\cd {\rm tr}G$, and let
$\E=\{\Gamma_1,\ldots,\Gamma_\ell\}$.  We need to show that the class of
$\E$ does not depend on the 
choice of the plinth $M$, the point $\omega\in\Omega$, or the
index $i\in\{1,\ldots,k\}$. First we verify that, given $M$, the class of $\E$ is
independent of $\omega$ and $i$. Let
$\omega_1,\ \omega_2\in\Omega$, and let
$\K_1=\K_{\omega_1}(\E)$, $\K_2=\K_{\omega_2}(\E)$ be the
corresponding Cartesian
systems for $M$. As $M$ is a
transitive subgroup of $G$, there is some $m\in M$ such that $\omega_1m=\omega_2$ and 
\cite[Lemma~2.2]{recog} implies that $\K_1^m=\K_2$. Hence, an element
$K$ of $\K_1$ involves a non-trivial, full strip $X$, if and only if
the corresponding element
$K^m$ of $\K_2$ involves the non-trivial full strip $X^m$. Thus the
existence of a non-trivial strip involved in an element of
a Cartesian system corresponding to $\E$ is
independent of the choice of $\omega$. 
Moreover, $\sigma_i(K)<T_i$ holds for some
$K\in\K_1$ if and only if $\sigma_i(K^m)<T_i$ holds. Thus the number
$|\mathcal F_i(\E,M,\omega)|$ is independent of the choice $\omega$,
and, by Theorem~\ref{stripth}(a), it is independent of $i$. 
It remains to prove that the definitions of the classes $\cd
{2\sim}G$ and $\cd{2\not\sim}G$ are independent of $\omega$ and of~$i$.
The simple argument given after~\eqref{f} shows that, for $j=1,\ 2$, the elements of
$\mathcal F_i(\E,M,\omega_j)$ are $G_{\omega_j}$-conjugate for some $i$ if and only if
they are $G_{\omega_j}$-conjugate for all $i$. Suppose that the
elements of $\mathcal F_1(\E,M,\omega_1)$ are $G_{\omega_1}$-conjugate. Let $\mathcal
F_1(\E,M,\omega_1)=\{A,B\}$ and let $g\in G_{\omega_1}$ be such that $A^g=B$. 
Then $\mathcal F_1(\E,M,\omega_1)^m=\mathcal F_1(\E,M,\omega_2)$, and $A^{mg^m}=A^{gm}=B^m$
hence the elements $A^m$ and $B^m$ of $\mathcal F_{1}(\E,M,\omega_2)$ are conjugate
under the element $g^m$. Further,
$\omega_2{g^m}=\omega_1{mg^m}=\omega_1{gm}=\omega_1m=\omega_2$, and so
$g^m\in G_{\omega_2}$. Thus the elements of $\mathcal F_{1}(\E,M,\omega_2)$ are
conjugate under $G_{\omega_2}$. Hence the definitions of $\cd {2\sim}G$
and $\cd{2\not\sim}G$ are independent of the choice of $\omega$ and~$i$. 

Next we show that the class of
$\E$ is independent of the plinth $M$. If $M$ is the unique
transitive, minimal
normal subgroup of $G$, then there is nothing to prove. If this is not
the case, then $G$ has exactly two transitive, minimal normal
subgroups $M$ and $\widehat M$, they are isomorphic, they both are 
regular on $\Omega$,
and $\cent G{M}=\widehat M$, $\cent G{\widehat M}=M$ also hold. Moreover, there
is an involution $\alpha\in\norm{\sym\Omega}G$  that interchanges $M$ and
$\widehat M$ (see~\cite[Lemma~5.1]{bp}). It follows from the proof
of~\cite[Lemma~5.1]{bp} that this involution $\alpha$ lies in a point
stabiliser. As $G$ is transitive on
$\Omega$, we may assume without loss of generality that $\alpha$ is an element
of the point stabiliser
$(\norm{\sym\Omega}G)_\omega$. Let
$\widehat T_1=T_1^\alpha,\ldots,\widehat T_k=T_k^\alpha$. Then
$\widehat M=\widehat T_1\times\cdots\times\widehat T_k$. 
Let $\widehat\sigma_i:\widehat M\rightarrow \widehat T_i$ be the natural
projection map, and define $\mathcal F_1(\E,M,\omega)$, 
$\mathcal F_1(\E,\widehat M,\omega)$ as in~\eqref{f}. Let $\K$ and
$\widehat \K$ be the Cartesian systems for $M$ and $\widehat M$,
respectively, with respect to $\omega$. For $i=1,\ldots,\ell$, let
$\gamma_i\in\Gamma_i$ such that $\omega\in\gamma_i$. As
$\omega\alpha=\omega$, we obtain that $\widehat
M_{\gamma_i}=(M_{\gamma_i})^\alpha$, and so, by the definition of $\K$
and $\widehat\K$, we have $\widehat \K=\K^\alpha=\{K^\alpha\
|\ K\in \K\}$. Let $t_1\cdots t_k\in M$ with $t_1\in T_1,\ldots,t_k\in
T_k$. Then $t_1^\alpha\cdots t_k^\alpha\in\widehat M$ with
$t_1^\alpha\in\widehat T_1,\ldots,t_k^\alpha\in\widehat T_k$. Thus
$\sigma_i(t_1\cdots t_k)^\alpha=\widehat\sigma_i((t_1\cdots
t_k)^\alpha)$ for
all $i\in\{1,\ldots,k\}$.
Therefore 
\begin{equation}\label{alpha}
\widehat\sigma_1(K^\alpha)=\sigma_1(K)^\alpha
\quad\mbox{for all}\quad K\in\K,
\end{equation} 
and so $\sigma_1(K)\neq T_1$ if and
only if $\widehat\sigma_1(K^\alpha)\neq \widehat T_1$. 
Hence $|\mathcal
F_1(\E,M,\omega)|=|\mathcal
F_1(\E,\widehat M,\omega)|$. It also follows that if $X$ is a
non-trivial, full strip involved in $K\in\K$, then $X^\alpha$ is a
non-trivial full strip involved in $K^\alpha\in\widehat \K$. Finally,
suppose that $\mathcal
F_1(\E,M,\omega)=\{A,B\}$ and $\mathcal
F_1(\E,\widehat M,\omega)=\{\widehat A,\widehat
B\}$. Equation~\eqref{alpha} implies that we may assume without loss of
generality that $\widehat A=A^\alpha$ and $\widehat B=B^\alpha$. 
Hence subgroups $A$ and $B$
are conjugate under an element $g\in G_\omega$, if and only if
$\widehat A$ ad $\widehat B$ are conjugate under $g^\alpha$, and,
since $g\in G_\omega$ and $\alpha\in(\norm{\sym\Omega}G)_\omega$, we
have $g^\alpha\in G_\omega$.  This
shows that the definitions of the classes $\cd 1G$, $\cd SG$, $\cd {1S}G$,
$\cd{2\sim}G$, $\cd {2\not\sim}G$, and $\cd 3G$ do not depend on the
choice of the plinth.

It is easy to see that the above classes are pairwise
disjoint, and Theorem~\ref{stripth} implies that each Cartesian
decomposition of $\cd {\rm tr}G$ belongs to at least one of the
classes. If $k=1$ and $M$ is a simple group,
then the map $\sigma_1$ is the identity map $M\rightarrow M$, and so
$\K_\omega(\E)=\mathcal F_1(\E,M,\omega)$. Theorem~6.1 of~\cite{recog}
implies that $|\E|=2$ and hence that $|\mathcal F_1(\E,M,\omega)|=2$. 
As $\E\in\cd{\rm tr}G$ we have that $G_\omega$
is transitive on $\K_\omega(\E)$, and hence on $\mathcal F_1(\E,M,\omega)$. Thus
$\E\in\cd{2\sim}G$. 
\end{proof}

The examples given earlier in this paper show that for each ${\rm
x}\in\{1,\ {\rm S},\ {\rm 1S},\ {2\!\sim},\ 2\!\not\sim,\ 3\}$ there
exists a group $G$ such that the
class $\cd xG$ is non-empty. Considering the corresponding
Cartesian systems, one can easily see that the the Cartesian
decomposition in the first example of the
introduction (given after Corollary~\ref{sdcor}) is in $\cd 1{T\wr K}$,
the decomposition of the second is in $\cd{2\sim}{(\aut{\alt 6})\wr K}$. 
The Cartesian decomposition in Example~\ref{ex2} belongs to
$\cd{2\not\sim}G$, the one in Example~\ref{ex3} is in $\cd 3G$, and
the one Example~\ref{ex1s} belongs to $\cd {1S}G$. Finally the decomposition
in the example given for
Theorem~\ref{diagonalembed} in the introduction is in $\cd S{K\wr
Q}$. On the other hand no group $G$ exists such that each of the classes $\cd
xG$ is non-empty. Indeed, it will follow from Theorems~\ref{last6}
and~\ref{diagonal} that
$\cd SG\neq\emptyset$ implies $\cd 1G\cup\cd{1S}G\cup\cd
{2\sim}G\cup\cd {2\not\sim}G\cup\cd 3G=\emptyset$, and $\cd 1G\cup\cd{1S}G\cup\cd
{2\sim}G\cup\cd {2\not\sim}G\cup\cd 3G\neq\emptyset$ implies  $\cd
SG=\emptyset$, for all  $G$.

If an observant reader
compares~\cite[Theorem~6.1]{recog} to our 6-Class Theorem, then
she finds that 
the Cartesian decompositions in~\cite[Theorem~6.1(ii)]{recog}
resemble the ones in $\cd 3G$. These Cartesian decompositions are,
however,  $G_\omega$-intransitive. This shows that allowing the
plinth $M$ to be non-simple leads to a very rich theory, and it would
not be realistic to expect results as explicit as the ones in~\cite{recog}.
Nevertheless, it is both 
possible and desirable to give further details of the Cartesian decompositions in
each class. 
The Cartesian decompositions in $\cd 3G$ were already described
in~\cite{design}. In this paper we describe classes $\cd 1G$ and $\cd
SG$. Besides the fact that these classes are, in a sense, the easiest,
their description also leads to the general result in
Theorem~\ref{diagonalembed}. Thus we felt that their description belongs in
this article. The classes $\cd {1S}G$, $\cd {2\sim}G$, and $\cd{2\not\sim}G$ pose considerably
more challenge, and they are addressed in a 
separate paper~\cite{1s}. 

Theorem~\ref{stripth} shows
that the existence  of a $G$-invariant Cartesian decomposition in a
particular class may pose severe restriction on the structure of
$G$. This is made more explicit in the following theorem. Recall that
an innately transitive group $G$ with plinth $M$ has
compound diagonal type if a point stabiliser 
$M_\omega$ is a subdirect subgroup of $M$
and is not simple.

\begin{theorem}\label{last6}
Suppose that $G$ is an innately transitive permutation group with a
non-abelian plinth $M$, and let $T$ be the isomorphism type of a
simple direct factor of $M$. Then the following all hold.
\begin{enumerate}
\item[(a)] If $\cd SG\neq\emptyset$ then $G$ is a quasiprimitive group
with compound diagonal type.
\item[(b)] If $\cd {1S}G\cup\cd {2\sim}G\neq\emptyset$ 
then $T$ admits a factorisation
with two proper, isomorphic subgroups, and hence $T$ is isomorphic to one of the
groups $\alt 6$, $\mat{12}$, $\pomegap 8q$, or $\sp 4{2^a}$ with $a\geq
2$. 
\item[(c)] If $\cd {2\not\sim}G\neq\emptyset$ then $T$ admits a
factorisation with proper subgroups.
\item[(d)] If $\cd 3G\neq\emptyset$ then $T$ admits a strong multiple
factorisation, and hence $T$ is isomorphic to one of the groups
$\sp{4a}2$ with $a\geq 2$, $\pomegap 83$, or $\sp 62$. 
\end{enumerate}
Moreover, for each ${\rm x}\in\{{\rm S},\ 1,\ {\rm 1S},\ 2\!\sim,\
2\!\not\sim,\ 3\}$ there is some $G$ as above such that $\cd
xG\neq\emptyset$. 
\end{theorem}
\begin{proof}
(a) Suppose now that $\cd{\rm S}G\neq\emptyset$, 
and let $\E\in\cd{\rm S}G$.  Then the elements of
$\K_\omega(\E)$ are proper subdirect subgroups of $M$. By
Theorem~\ref{stripth}, 
the set $\mathcal P$ of supports $\supp X$  of non-trivial strips $X$
involved in
$\K_\omega(\E)$ is non-empty and, by Theorem~\ref{stripth}(d), $\mathcal P$
is a $G_\omega$-invariant partition of
$\{T_1,\ldots,T_k\}$. 
Hence each $K\in\K_\omega(\E)$ is a direct product of those
non-trivial strips that are involved in $K$, and those of the $T_i$
that are not involved in these strips. Therefore  
$$
M_\omega=\bigcap_{K\in\K_\omega(\E)}K=\prod_{\hbox{\scriptsize$X$ is a nontrivial strip
involved in $\K_\omega(\E)$}}X.
$$
Thus $M_\omega$ is a subdirect subgroup of $M$, and hence $G$ has
diagonal type. It follows from~\cite[Proposition~5.5]{bp} that $G$ is
quasiprimitive. Since there are at
least $|\K_\omega(\E)|\geq 2$ strips $X$ involved in $\K_\omega(\E)$,
we have that $M_\omega\not\cong T$; therefore $G$ is of compound
diagonal type. 

(b)--(d) It follows from Theorem~\ref{stripth} and the definitions of
$\cd{2\sim}G$ and $\cd{2\not\sim}G$ that the claimed factorisations are admitted by
$T$. Factorisations of finite simple groups with isomorphic subgroups
were listed in~\cite[Lemma~5.2]{recog}, while strong multiple
factorisation of finite simple groups were classified
in~\cite[Table~V]{bad:fact}. These results imply that the isomorphism
type of $T$ in parts~(b) and~(d) is as claimed.

The final assertion follows from the discussion between
Theorems~\ref{5class} and~\ref{last6}
\end{proof}

There are finite simple groups that admit no
factorisations with proper subgroups, for example, 
${\sf PSU}_{2m+1}(q)$ (with a finite number of exceptions) and some 
sporadic groups (see the tables in~\cite{lps:max}). Therefore Theorem~\ref{last6}(c) also restricts the
isomorphism type of $T$. 

\section{The Cartesian decompositions in $\cd {\rm S}G$ and in $\cd 1G$}\label{5}

In
this section we 
describe 
the elements of $\cd {\rm S}G$ and $\cd {1}G$ for an innately transitive
permutation group $G$.
First we show how to construct such Cartesian decompositions. 
Theorem~\ref{c1}
implies that our construction is as
general as possible.

\begin{example}\label{exc1}
Suppose that $G$ is an innately transitive permutation group with
a non-abelian plinth $M$, and let $T_1,\ldots,T_k$ be the simple
normal subgroups of $M$.
Note that the conjugation action of $G$ permutes the set
$\{T_1,\ldots,T_k\}$ transitively.
Suppose that
$\omega\in\Omega$ and 
$\mathcal A=\{A_1,\ldots,A_\ell\}$ is a $G$-invariant partition of $\{T_1,\ldots,T_k\}$ such
that 
$$
M_\omega=\sigma_{A_1}(M_\omega)\times\cdots\times\sigma_{A_\ell}(M_\omega).
$$
Then for $i=1,\ldots,\ell$ let
$$
K_i=\sigma_{A_i}(M_\omega)\times\prod_{T_j\not\in A_i}T_j
$$
and set $\K=\{K_1,\ldots,K_\ell\}$. 

If $g\in G_\omega$ and $i\in\{1,\ldots,\ell\}$ then, since $M_\omega$
is normalised by $G_\omega$,
$$
\left(\sigma_{A_i}(M_\omega)\right)^g=\sigma_{A_i^g}(M_\omega^g)=\sigma_{A_i^g}(M_\omega).
$$
Therefore
$$
K_i^g=\left(\sigma_{A_i}(M_\omega)\times\prod_{T_j\not\in
A_i}T_j\right)^g=
\sigma_{A_i^g}(M_\omega)\times\prod_{T_j\not\in A_i^g}T_j,
$$
and so $K_i^g\in\K$. Hence $\K$ is
$G_\omega$-invariant. It is also easy to see that
equations~\eqref{csdef1} and~\eqref{csdef2} hold for
$\K$, and hence $\K$ is a Cartesian system for $M$ with respect to
$\omega$.  It follows from the last displayed equation that the
$G_\omega$-actions on $\K$ and on $\mathcal A$ are equivalent.  Hence
$G_\omega$ is transitive on $\K$, and so $\E(\K)\in\cd{\rm tr}G$. 
Moreover $\E(\K)\in\cd{\rm S}G$ if $M_\omega$ is a
subdirect subgroup of $M$, and $\E(\K)\in\cd{1}G$ otherwise.
\end{example}

\begin{theorem}\label{c1}
Let $G$ be an innately transitive group on $\Omega$ with a non-abelian
plinth $M$, and let $\omega\in\Omega$. Let ${\rm x}=\rm S$ if $M_\omega$ is a subdirect subgroup of $M$, and let
${\rm x}=1$ otherwise. Then there is a bijection between the set 
$\cd xG$
and the set 
of $G$-invariant partitions $\{P_1,\ldots,P_\ell\}$ of
$\{T_1,\ldots,T_k\}$ satisfying
$M_\omega=\sigma_{P_1}(M_\omega)\times\cdots\times\sigma_{P_\ell}(M_\omega)$,
and each $\E\in\cd xG$ arises as in Example~$\ref{exc1}$. 
\end{theorem}
\begin{proof}
If $\mathcal P=\{P_1,\ldots,P_\ell\}$ is a $G$-invariant partition of
$\{T_1,\ldots,T_k\}$ such that
$$
M_\omega=\sigma_{P_1}(M_\omega)\times\cdots\times\sigma_{P_\ell}(M_\omega),
$$
then let $\K(\mathcal P)=\{K_1,\ldots,K_\ell\}$ where, as in Example~\ref{exc1},
$$
K_i=\sigma_{P_i}(M_\omega)\times\prod_{T_j\not\in
P_i}T_j\quad\mbox{for all}\quad i\in\{1,\ldots,\ell\}.
$$
Then, arguing as in Example~\ref{exc1}, $\K(\mathcal P)$ is a
$G_\omega$-invariant Cartesian system 
of subgroups for $M$. It is easy to check that the map $\mathcal P\mapsto\K(\mathcal
P)$ is injective. If $M_\omega$ is not a subdirect subgroup of $M$,
then $\E(\K(\mathcal P))\in\cd 1G$, and $\E(\K(\mathcal P))\in\cd{\rm
S}G$ otherwise. Hence we have an injective map $\mathcal
P\mapsto\E(\K(\mathcal P))$ from the set of $G$-invariant partitions of
$\{T_1,\ldots,T_k\}$ to $\cd xG$. 

Now let $\E\in\cd{1}G$ and set $\K_\omega(\E)=\{K_1,\ldots,K_\ell\}$ and 
$\mathcal P=\{\{T_i\ |\ \sigma_i(K_j)\neq T_i\}\ |\
j=1,\ldots,\ell\}$. Since $\E\in\cd 1G$ we have $|\mathcal F_i(\E,M,\omega)|=1$ for
all $i$, that is there is a unique $j$ such that $\sigma_i(K_j)<T_i$. Thus $\mathcal P$
is a partition of $\{T_1,\ldots,T_k\}$. As $\E\in\cd 1G$, it follows
from Lemma~\ref{scottgen} that if $i\in\{1,\ldots,k\}$ and
$j\in\{1,\ldots,\ell\}$ such that $\sigma_i(K_j)=T_i$, then $T_i\leq
K_j$. This implies that
$\K_\omega(\E)=\K(\mathcal P)$, and we have $\E=\E(\K(\mathcal P))$. 
Thus in this case $\mathcal P\mapsto\E(\K(\mathcal
P))$ is surjective, 
and so it is a bijection. On the other hand, if $M_\omega$ is a subdirect subgroup of
$M$ and $\E\in\cd{\rm S}G$, then Theorem~\ref{disstrips} shows that
$\mathcal P=\{\{T_i\ |\ T_i\not\leq K_j\}\ |\ j=1,\ldots,\ell\}$
is a partition of $\{T_1,\ldots,T_k\}$.
It follows from the definition of $\mathcal P$ that
$\K_\omega(\E)=\K(\mathcal P)$, and hence again $\E=\E(\K(\mathcal
P))$, and   the map $\mathcal P\mapsto\E(\K(\mathcal
P))$ is a bijection in this case also.
\end{proof}

Using the above characterisation of Cartesian decompositions, we 
can describe the class of Cartesian decompositions
that are preserved by innately transitive groups with diagonal
type, which, we recall, are defined as follows. Suppose that
$G$ is an innately transitive permutation group on $\Omega$ with a non-abelian
and non-simple plinth $M=T_1\times\cdots\times T_k$, where the $T_i$ are
isomorphic to a non-abelian simple group $T$.  As above, let $\sigma_i$ denote the $i$-th coordinate
projection $M\rightarrow T_i$. Let $\omega\in\Omega$. Then $G$ has diagonal type if
$\sigma_i(M_\omega)=T_i$ for some (and hence all) $i\in\{1,\ldots,k\}$. 
It follows from Scott's lemma that, in this case, $M_\omega$ 
is isomorphic to $T^m$ for some $m\leq k$. 
By~\cite[Proposition~5.5]{bp}, $M$ is the unique minimal normal
subgroup of $G$, and so $G$ is quasiprimitive.
We say that such a $G$ has {\em simple diagonal type} if
$M_\omega\cong T$, and {\em compound diagonal type} otherwise (see also~\cite{prae:quasi}).

\begin{theorem}\label{diagonal}
Let $G$ be an innately transitive permutation group of diagonal
type. Then $\cd {}G\neq \emptyset$ if and only if $G$ is a quasiprimitive
group with compound diagonal type. Moreover, if $G$ is of compound
diagonal type then $\cd {}G=\cd SG\neq \emptyset$. 
\end{theorem}
\begin{proof}
Let $M$ be the unique minimal normal subgroup of $G$. 
Let $T_1,\ldots,T_k$ be the simple normal subgroups of $M$, and let
$T$ denote their common isomorphism type. Suppose
that $M_\omega$ is a subdirect subgroup of $M$, and let
$\E$ be an element in $\cd{}G$. It follows
from~\cite[Proposition~5.5]{bp} that $G$ is a quasiprimitive group. 
For each $K\in\K_\omega(\E)$ we have $M_\omega\leq K$, and
so 
all elements of $\K_\omega(\E)$ are subdirect subgroups of $M$. Let
$K_1,\ K_2\in\K_\omega(\E)$ be two distinct subgroups. Then $K_1,\
K_2\neq M$, and so, by Lemma~\ref{scott}, $K_1,\ K_2$ involve non-trivial full strips $X_1$ and
$X_2$, say, and, by Theorem~\ref{disstrips}, $X_1\neq X_2$.  Suppose that $T_{i_1}\in\supp X_1$ and $T_{i_2}\in\supp
X_2$. Then $T_{i_1}^g=T_{i_2}$ for some $g\in G_\omega$, and so
Theorem~\ref{disstrips} implies that $X_1^g=X_2$. This, in turn, yields
$K_1^g=K_2$. Thus $\E\in\cd{\rm tr}G$, and
then clearly $\E\in\cd{\rm S}G$. Hence $\cd{}G=\cd{\rm S}G$.  Note that our argument also implies
that if a non-trivial, full strip $X$ is involved in $\K$, then it is
also involved in $M_\omega$. 
The supports of
the full strips involved in $M_\omega$ form a $G$-invariant partition
$\mathcal P$ of $\{T_1,\ldots,T_k\}$ such that
$M_\omega=\prod_{P\in\mathcal P}\sigma_P(M_\omega)$. Since $|\mathcal
P|\geq 2$, $M_\omega$ is not simple, and so $G$ has compound diagonal
type. Conversely, if
$G$ has compound diagonal type, then, by Theorem~\ref{c1}, $\cd{}G$
is non-empty, and we showed above that, in this case, $\cd {}G=\cd
SG$. 
\end{proof}

\section{Inclusions in wreath products}\label{final}

This section is devoted to proving
Theorems~\ref{main},~\ref{diagonalembed} and Corollary~\ref{sdcor}. 
Throughout the 
section we assume that the common hypotheses of Theorems~\ref{main}
and~\ref{diagonalembed}
hold. 
Thus, let $H$ be a quasiprimitive,  almost simple permutation group 
acting on a set $\Gamma$, and let $U$ denote its simple normal subgroup.
Set $W=H\wr\sy \ell\cong H^\ell\rtimes \sy\ell$, and consider $W$ as a permutation
group on $\Gamma^\ell$ in product action. 
Let $N=U_1\times\cdots\times U_\ell\cong U^\ell$ denote the unique
minimal normal subgroup of
$W$. As there is a natural isomorphism $U\rightarrow U_i$ for each $i\in\{1,\ldots,\ell\}$, we
consider the $U_i$ as subgroups of $\sym\Gamma$. Let $G$ be a finite,  innately transitive permutation group
acting on $\Gamma^\ell$ with
a non-abelian plinth $M=T_1\times\cdots\times T_k$ where $T_1,\ldots,T_k$ are
finite simple groups all isomorphic to a group $T$. 
Assume that $G\leq W$. Note that at this stage we do not assume that
$G$ projects onto a transitive subgroup of $\sy\ell$. 

Let us introduce some extra notation to facilitate our investigation
of the
situations described by Theorems~\ref{main} and~\ref{diagonalembed}. As before, $\sigma_i$
denotes the natural projection $M\rightarrow T_i$, and, for
$i=1,\ldots,\ell$, 
let $\mu_i$
denote the natural projection $N\rightarrow U_i$. Note that $W$
preserves the natural Cartesian decomposition
$\E=\{\Gamma_1,\ldots,\Gamma_\ell\}$
of $\Gamma^\ell$
where 
$$
\Gamma_i=\{\{(\gamma_1,\ldots,\gamma_\ell)\ |\ \gamma_i=\gamma\}\ |\
 \gamma\in\Gamma\}\quad\mbox{for}\quad i\in\{1,\ldots,\ell\}.
$$ 
Fix
$\gamma\in\Gamma$ and let $\omega=(\gamma,\ldots,\gamma)$. Let $\K=\{\K_1,\ldots,K_\ell\}$
be the $G_\omega$-invariant Cartesian system in $M$ with respect
to $\omega$ corresponding to $\E$. It follows from the definition of
the product action that
the $W$-actions on $\E$ and on the set of simple direct factors
$U_1,\ldots,U_\ell$ of $N$ are equivalent.

The proof of the next lemma uses the Schreier Conjecture that the
outer automorphism group of a finite simple group is soluble. The
validity of the Schreier Conjecture
follows from the finite simple group classification.

\begin{lemma}\label{minn}
We have that $M\leq N$.
\end{lemma}
\begin{proof}
Let $B$ denote the base group $H^\ell$ of $W$. It follows from the
definition of the product action that the pointwise stabiliser in $W$ of
$\E$ coincides with $B$, and
so~\cite[Proposition~2.1]{recog} implies that $M\leq B$. Now, as $M$ is a
minimal normal subgroup of $G$ we have that either $M\leq N$ or $M\cap
N=1$. If $M\leq N$, then we are done, so assume that
$M\cap N=1$. Then 
$(MN)/N\cong M/(M\cap N)\cong M$, and so $M$ can be viewed as a subgroup
of $B/N$. On the other hand, $B/N\cong (H/U)^\ell$. By the Schreier
Conjecture, $H/U$ is a soluble group, and therefore so is $B/N$. Hence assuming
that $M\cap N=1$ leads to the incorrect conclusion that $M$ is a soluble
group. Therefore $M\leq N$ must hold.
\end{proof}

Note that $M$ is a characteristically simple group and each of its normal
subgroups is a product of some of the $T_i$. If $L$ is a normal
subgroup of $M$ then the quotient $M/L$ can naturally be identified
with the subgroup $\prod_{T_i\not\leq L}T_i$, and in future 
this identification
will be used without further comment. For a subgroup $K<M$, ${\sf
Core}_M(K):=\bigcap_{m\in M}K^m$ is the largest normal subgroup of $M$
contained in $K$.

\begin{lemma}\label{mu}
Let $i\in\{1,\ldots,\ell\}$. Then 
\begin{equation}\label{muM}
M^{\Gamma_i}= \prod_{T_j\not\leq K_i}T_j\cong T^{s_i}
\end{equation}
for some $s_i\geq 1$, and
$\mu_i(M)$ is permutationally
isomorphic to $M^{\Gamma_i}$. Moreover, $\mu_i(M)$ is a
transitive subgroup of $U_i$, and if $s_i\geq 2$ then $U=\Alt\Gamma$.
\end{lemma}
\begin{proof}
Note that $K_i$
is the stabiliser of a point for the $M$-action on $\Gamma_i$, and the kernel of this action is ${\sf Core}_{M}(K_i)$. Now $T_j\leq {\sf
Core}_{M}(K_i)$ if and only if $T_j\leq K_i$, and so~\eqref{muM}
holds.

Let $\alpha$ denote the bijection
$\Gamma\rightarrow\Gamma_i$ mapping
$$
\gamma\mapsto\{(\gamma_1,\ldots,\gamma_\ell)\ |\
\gamma_i=\gamma\}\quad\mbox{for all}\quad\gamma\in\Gamma.
$$
The map $\mu_i$ can be considered as a permutation representation of
$M$ in $\sym\Gamma$. We claim that the $M$-actions on $\Gamma$ and
$\Gamma_i$ are equivalent via the bijection $\alpha$. Let
$\gamma\in\Gamma$ and $m\in M$. Then
$$
\alpha(\gamma m)=\{(\gamma_1,\ldots,\gamma_\ell)\ |\
\gamma_i=\gamma m\}=
\{(\gamma_1,\ldots,\gamma_\ell)\ |\
\gamma_i=\gamma\}m=\alpha(\gamma)m.
$$ 
Therefore our claim holds. As $\Gamma_i$ is an $M$-invariant partition
of $\Omega$, the group $M$ is transitive on $\Gamma_i$, and so $\mu_i(M)$ is
transitive on $\Gamma$. Thus the factorisation
$\mu_i(M)(U_i)_\gamma=U_i$ holds. Also note that, $\mu_i(M)$ is a
homomorphic image of $M$, and so $\mu_i(M)\cong T^{s_i}$ for some
$s_i$. Let $I$ denote the set of indices $i$ such that $s_i\geq 1$. 
Then $M\leq\prod_{i\in I}U_i$, and so $\prod_{i\in I}U_i$ is a
transitive normal subgroup of $N$. However, by the definition of the
product action, $\prod_{j\in J}U_j$ is intransitive if $J$ is a proper
subset of $\{1,\ldots,\ell\}$. Hence $I=\{1,\ldots,\ell\}$, and so
$s_i\geq 1$ for all~$i$.
If $s_i\geq 2$ for some~$i$, 
then~\cite[Theorem~1.4]{bad:quasi} implies that $U=\Alt\Gamma$. 
\end{proof}

These results enable us to show that, for four of the six classes
identified by the 6-Class Theorem, the corresponding embedding of $G$
belongs to Theorem~\ref{main}(c). 

\begin{lemma}\label{scasec}
Suppose that $\E\in\cd{tr}G$. If $\mu_i(M)\cong T^{s_i}$ with $s_i\geq
2$ for some
$i\in\{1,\ldots,\ell\}$, then Theorem~$\ref{main}$(c)
is valid.
\end{lemma}
\begin{proof}
Suppose that $\mu_i(M)\cong T^{s_i}$ and $s_i\geq 2$ for some $i\geq 1$. It follows from Lemma~\ref{mu} that
in this case $U=\Alt\Gamma$. Also, note that there are exactly $s_i$ indices $j$
such that $\mu_i(T_j)\cong T$. 
On the other hand, in cases~(a) and (b) for
each $i\in\{1,\ldots,\ell\}$ there is a unique $j\in\{1,\ldots,k\}$
such that $\mu_i(T_j)\cong T$. Therefore Theorem~\ref{main}(c) holds.
\end{proof}

\begin{lemma}\label{dcor}
(a) If $\E\in \cd SG\cup\cd {1S}G\cup\cd 3G\cup\cd {2\not\sim}G$,
then $\mu_i(M)\cong T^s$ with $s$ independent of $i$ and $s\geq 2$,
and Theorem~$\ref{main}(c)$ is valid.

(b) If $\E\in\cd 1G$ then either case~(a) or case~(c) of
Theorem~$\ref{main}$ is valid.

(c) If $\E\in\cd {2\sim}G$ then either case~(b) or case~(c) of
Theorem~$\ref{main}$ is valid.
\end{lemma}
\begin{proof}
In each case, $G$ is transitive
on $\E$, and $M$ is a normal subgroup of $G$, and so the permutation
groups $M^{\Gamma_i}$ are pairwise permutationally 
isomorphic. Thus $M^{\Gamma_i}\cong
T^s$ for some $s\geq 1$ independent of $i$, and, by Lemma~\ref{mu},
$\mu_i(M)\cong T^s$. 
If $s\geq 2$ then
Lemma~\ref{scasec} implies 
that case~(c) of the theorem holds. 

(a) Let $\E\in \cd SG\cup\cd {1S}G\cup\cd 3G\cup\cd {2\not\sim}G$. 
By the discussion in the previous paragraph, 
it suffices to verify that
$\mu_1(M)\cong T^s$ with $s\geq 2$.
Note that, by Lemma~\ref{mu},
$\mu_1(M)=\prod_{T_j\not\leq K_1}T_j$. If $\E\in\cd SG\cup\cd
{1S}G$, then $K_1$ involves a non-trivial, full strip, and so $s\geq
2$ follows immediately. 

Suppose next that
$\E\in\cd 3G$. Then $\mathcal F_1(\E,M,\omega)=\{A,B,C\}$ for some subgroups $A$,
$B$, and $C$ of $T_1$, such that $A$, $B$, and $C$ form a strong
multiple factorisation of $T_1$. Hence~\cite[Table~V]{bad:fact} yields
that the subgroups $A$, $B$, and $C$ represent three
distinct isomorphism types. Then there are pairwise distinct indices
$j_1,\ j_2,\ j_3\in\{1,\ldots,\ell\}$ such that
$\sigma_1(K_{j_1})=A$, $\sigma_1(K_{j_2})=B$, $\sigma_1(K_{j_3})=C$. 
Let $g_1,\ g_2,\ g_3\in G_\omega$ such that
$K_{j_1}^{g_1}=K_{j_2}^{g_2}=K_{j_3}^{g_3}=K_1$. Let
 $i_1,\ i_2,\ i_3\in\{1,\ldots,k\}$ such that
$T_1^{g_1}=T_{i_1}$, $T_1^{g_2}=T_{i_2}$, $T_1^{g_3}=T_{i_3}$. Then
$\sigma_{i_1}(K_1)=\sigma_1(K_{j_1})^{g_1}=A^{g_1}$,
$\sigma_{i_2}(K_1)=\sigma_1(K_{j_2})^{g_2}=B^{g_2}$, and
$\sigma_{i_3}(K_1)=\sigma_1(K_{j_3})^{g_3}=C^{g_3}$. As $A$, $B$, and
$C$ are pairwise non-isomorphic, it follows that $i_1,\ i_2,\ i_3$ are
also pairwise distinct. Thus~\eqref{muM} implies that
$\mu_1(M)\cong T^s$ for some $s\geq
3$.

Suppose finally that $\E\in\cd {2\not\sim}G$ such that $\mathcal F_1(\E,M,\omega)=\{A,B\}$ where
$A,\ B\leq T_1$ are not conjugate under $G_\omega$. Then it follows
that there are indices $j_1,\ j_2\in\{1,\ldots,\ell\}$ such that
$\sigma_{1}(K_{j_1})=A$ and $\sigma_1(K_{j_2})=B$. Let $g_1,\ g_2\in
G_\omega$ be such that $K_{j_1}^{g_1}=K_{j_2}^{g_2}=K_1$. Then
$T_1^{g_1}=T_{i_1}$ and $T_1^{g_2}=T_{i_2}$ for some $i_1,\
i_2\in\{1,\ldots,k\}$. Hence
$\sigma_{i_1}(K_1)=\sigma_1(K_{j_1})^{g_1}=A^{g_1}$ and
$\sigma_{i_2}(K_1)=\sigma_1(K_{j_2})^{g_2}=B^{g_2}$. As $A$ and $B$
are not conjugate under $G_\omega$, we have that $A^{g_1}\neq
B^{g_2}$, and so $i_1\neq i_2$. 
Thus~\eqref{muM} implies that
$\mu_i(M)\cong T^s$ where $s\geq
2$, and the result follows.

(b) As noted above, if $\mu_i(M)\cong T^s$ with $s\geq 2$ then
Theorem~\ref{main}(c) holds. Assume now that
$\mu_i(M)\cong T$.
As $\E\in\cd 1G$,
for all $i\in\{1,\ldots,k\}$,
there is a unique $j\in\{1,\ldots,\ell\}$ such that $T_i\not\leq K_j$. This means that for
all $i\in\{1,\ldots,k\}$ there is a unique $j\in\{1,\ldots,\ell\}$ such that $\mu_j(T_i)\neq 1$, and so
$T_i\leq U_j$. On the other hand, as $\mu_i(M)\cong T$, for all $j\in\{1,\ldots,\ell\}$
there is a unique $i\in\{1,\ldots,k\}$ such that $T_i\leq U_j$.
Therefore
$\ell=k$ and the $T_i$ and the $U_j$ can be indexed
so that $T_1\leq U_1,\ldots,T_k\leq U_k$, and so Theorem~\ref{main}(a) holds.

(c) As in part~(b), either Theorem~\ref{main}(c) holds or
$\mu_i(M)\cong T$ for all $i$; assume the latter.
As $\E\in\cd{2\sim}G$, for all $i\in\{1,\ldots,k\}$, there are exactly two 
indices $j\in\{1,\ldots,\ell\}$ such that $\mu_j(T_i)\neq 1$. If
$j_1,\ j_2\in\{1,\ldots,\ell\}$ are these indices then we obtain that
$T_i\leq U_{j_1}\times U_{j_2}$. On the other hand, as $\mu_j(M)\cong T$, for all $j\in\{1,\ldots,\ell\}$
there is a unique $i\in\{1,\ldots,k\}$ such that $\mu_j(T_i)\neq
1$. 
Counting the pairs in the set 
$$
\{(i,j)\ |\ i\in\{1,\ldots,k\},\ j\in\{1,\ldots,\ell\},\ T_i\not\leq
K_j\}
$$
we obtain that $\ell=2k$ and the $T_i$ and the $U_i$ can be indexed
so that $$
T_1\leq U_1\times U_2,\ T_2\leq U_3\times U_4,\ldots,T_k\leq
U_{2k-1}\times U_{2k}. 
$$

Note that $U_1\times U_2$ can be viewed as a permutation group acting on
$\Gamma_1\times \Gamma_2$ preserving the  natural Cartesian
decomposition formed by the $(U_1\times U_2)$-invariant partitions
$\Gamma_1$ and $\Gamma_2$. Choose $\gamma_1\in\Gamma_1$,
$\gamma_2\in\Gamma_2$ and set $\omega=(\gamma_1,\gamma_2)$. 
Then~\cite[Theorem~6.1]{recog} implies that
the isomorphism type $T$ of $T_1$, and the stabilisers $(T_1)_{\gamma_1}$,
$(T_1)_{\gamma_2}$ and $(T_1)_{\omega}$ are as in
Table~\ref{table2}. Hence $T_1$ acts
primitively on both $\Gamma_1$ and $\Gamma_2$. Thus $(U_1)^{\Gamma_1}$ and $(U_2)^{\Gamma_2}$
are primitive permutation groups. 
Since $(T_1)^{\Gamma_1}\leq U_1^{\Gamma_1}$ and $(T_2)^{\Gamma_2}\leq U_2^{\Gamma_2}$, 
the results of~\cite{lps:maxsub}
yield that $U$ and
 $T$ are as in the corresponding columns of
Table~\ref{maintable}. Thus Theorem~\ref{main}(b) holds.
\end{proof}
\begin{center}
\begin{table}
$$
\begin{array}{|c|c|c|c|}
\hline
& T & (T_1)_{\gamma_1},\ (T_1)_{\gamma_2}& (T_1)_{\omega}\\
\hline\hline
1& \alt 6 & \alt 5,\ \tau(\alt 5)\ (\tau\not\in\sy 6) & \dih{10}\\
\hline
2 & \mat{12} & \mat{11},\ \tau(\mat{11})\ (\tau\not\in\mat{12}) & \psl 2{11} \\
\hline
3 & \pomegap 8q &  \Omega_7(q),\ \tau(\Omega_7(q))\mbox{
($\tau$ a triality)} & {\sf G}_2(q)\\
\hline
4 & \sp 4{2^a},\ a\geq 2 & \sp{2}{2^{2a}}\cdot 2,\ {\sf O}^-_4(q)
&\dih{q^2+1}\cdot 2\\
\hline
\end{array}
$$
\caption{The stabilisers in $T_1$}\label{table2}
\end{table}
\end{center}

The proof of Theorem~\ref{main} is now very easy.

\begin{proof}[Proof of Theorem~$\ref{main}$]
The theorem follows from the 6-Class Theorem (Theorem~\ref{5class})
and Lemma~\ref{dcor}.
\end{proof}

Next we prove Theorem~\ref{diagonalembed}.

\bigskip

\begin{proof}[Proof of Theorem~$\ref{diagonalembed}$]
Let $G$ and $W$ be as in Theorem~\ref{diagonalembed}. Then
$\cd{}G\neq\emptyset$, and hence, by Theorem~\ref{diagonal}, $G$ is a
quasiprimitive group with compound diagonal type and the Cartesian
decomposition $\E$ corresponding to the product action of $W$ belongs
to $\cd{S}G$. In particular, $\E\in\cd{tr}G$, and so $G$ projects onto
a transitive subgroup of $\sy\ell$. Thus part~(a) is proved. It
follows from Lemma~\ref{dcor} that the inclusion $G\leq W$ is as in
Theorem~\ref{main}(c), and so $U=\Alt\Gamma$. Thus part~(b) holds.

By Lemma~\ref{minn}, $M\leq N$.
Let $\omega\in\Omega$ and consider the Cartesian system
$\K_\omega(\E)$ for $M$. 
If $K\in\K_\omega(\E)$ then $M_\omega\leq K$, and so $K$ is a
subdirect subgroup of $M$. Thus Lemma~\ref{scott} implies that $K$ is
the direct product of pairwise disjoint, full strips. 
By Theorem~\ref{disstrips}, if
$X_1$ and $X_2$ are non-trivial, full strips involved in $\K_\omega(\E)$ then
$X_1$ and $X_2$ are disjoint. Thus if
$T_i\not\leq K_j$ for some $i\in\{1,\ldots,k\}$ and
$j\in\{1,\ldots,\ell\}$ then $T_i\leq K_{j'}$ for all
$j'\in\{1,\ldots,\ell\}\setminus\{j\}$. Hence the
$T_1,\ldots,T_k$ and the $U_1,\ldots,U_\ell$ can be indexed such that 
$$
M^{\Gamma_1}=T_1\times\cdots\times
T_m,\ldots,M^{\Gamma_\ell}=T_{(\ell-1)m+1},\ldots,T_{\ell m}.
$$
Thus $k=\ell m$ and by Lemmas~\ref{mu} and~\ref{dcor}, $m\geq 2$. Also
$$
T_1\times\cdots\times
T_m\leq U_1,\ldots,T_{(\ell-1)m+1}\times\ldots\times T_{\ell m}\leq U_\ell,
$$ 
and statement~(c) is also valid.
\end{proof}

Finally we prove Corollary~\ref{sdcor}.

\begin{proof}[Proof of Corollary~$\ref{sdcor}$]
Let $G$ be an innately transitive group on $\Omega$ 
with simple diagonal type, and let $M$ be the unique minimal normal subgroup
of $G$, whose existence is guaranteed by~\cite[Proposition~5.5]{bp}. Assume
that the result is not true for $G$, so $G\leq W$ where $W$ is permutationally
isomorphic to $\sym\Gamma\wr\sy\ell$ with some $|\Gamma|\geq 2$ and
$\ell\geq 2$. Let $\E$ denote the natural Cartesian decomposition
of $\Omega$ corresponding to the product action of
$W$. Then~\cite[Proposition~2.1]{recog} implies that $M$ lies in the
pointwise stabiliser in $W$ of $\E$, and so $M$ is a subgroup of the
base group $(\sym\Gamma)^\ell$ of $W$. As $M$ is a non-abelian
characteristically simple
group and $\sym\Gamma$ is soluble for $|\Gamma|\leq 4$, we obtain that $|\Gamma|\geq 5$. Thus
Theorem~\ref{diagonalembed} implies that $G$ has compound diagonal
type, which is a contradiction. 
\end{proof}

\end{document}